\documentclass[12pt]{amsart}
\usepackage{a4wide}
\usepackage{amssymb}
\usepackage{amsthm}
\usepackage{amsmath}
\usepackage{amscd}
\usepackage[mathscr]{eucal}
\usepackage[all]{xy}

\numberwithin{equation}{section}

\theoremstyle{plain}
\newtheorem{theorem}{Theorem}[section]
\newtheorem{corollary}[theorem]{Corollary}

\newtheorem{proposition}[theorem]{Proposition}
\newtheorem{assumption}[theorem]{Assumption}

\theoremstyle{definition}

\theoremstyle{remark}

\newcommand{\R}{\mathbb{R}}
\newcommand{\Q}{\mathbb{Q}}
\newcommand{\Z}{\mathbb{Z}}
\newcommand{\N}{\mathbb{N}}
\newcommand{\C}{\mathbb{C}}

\renewcommand{\H}{\mathbb{H}}

\newcommand{\D}{\mathbf{D}}

\newcommand{\zxz}[4]{\begin{pmatrix} #1 & #2 \\ #3 & #4 \end{pmatrix}}
\newcommand{\abcd}{\zxz{a}{b}{c}{d}}

\newcommand{\kzxz}[4]{\left(\begin{smallmatrix} #1 & #2 \\ #3 & #4\end{smallmatrix}\right) }
\newcommand{\kabcd}{\kzxz{a}{b}{c}{d}}

\newcommand{\calD}{\mathcal{D}}

\newcommand{\calF}{\mathcal{F}}

\newcommand{\calH}{\mathcal{H}}

\newcommand{\calL}{\mathcal{L}}
\newcommand{\calM}{\mathcal{M}}
\newcommand{\calN}{\mathcal{N}}
\newcommand{\calO}{\mathcal{O}}

\newcommand{\calT}{\mathcal{T}}

\newcommand{\calW}{\mathcal{W}}
\newcommand{\calX}{\mathcal{X}}

\newcommand{\frakd}{\mathfrak d}

\newcommand{\frakg}{\mathfrak g}

\newcommand{\eps}{\varepsilon}
\newcommand{\bs}{\backslash}
\newcommand{\norm}{\operatorname{N}}
\newcommand{\vol}{\operatorname{vol}}

\newcommand{\Sl}{\operatorname{SL}}

\newcommand{\Symp}{\operatorname{Sp}}

\newcommand{\Orth}{\operatorname{O}}
\newcommand{\Uni}{\operatorname{U}}
\newcommand{\Hom}{\operatorname{Hom}}
\newcommand{\Spec}{\operatorname{Spec}}
\newcommand{\Proj}{\operatorname{Proj}}
\newcommand{\End}{\operatorname{End}}
\newcommand{\sym}{\text{\rm sym}}
\newcommand{\pre}{\text{\rm pre}}
\newcommand{\Pet}{\text{\rm Pet}}

\newcommand{\Gal}{\operatorname{Gal}}

\newcommand{\supp}{\operatorname{supp}}

\newcommand{\Pic}{\operatorname{Pic}}

\newcommand{\ord}{\operatorname{ord}}

\newcommand{\dv}{\operatorname{div}}
\newcommand{\OK}{\mathcal{O}_K}
\newcommand{\GK}{\Gamma_K}
\newcommand{\ch}{\operatorname{CH}}
\newcommand{\z}{\operatorname{Z}}
\newcommand{\za}{\operatorname{\widehat{Z}}}
\newcommand{\rata}{\operatorname{\widehat{Rat}}}
\newcommand{\dega}{\operatorname{\widehat{deg}}}
\newcommand{\pica}{\operatorname{\widehat{Pic}}}

\newcommand{\cha}{\operatorname{\widehat{CH}}}
\newcommand{\cc}{\operatorname{c}}
\newcommand{\cca}{\widehat{\operatorname{c}}}
\newcommand{\fh}{\operatorname{ht}}

\begin{document}

\title{Infinite products in number theory and geometry}
\author[Jan H.~Bruinier]{Jan Hendrik Bruinier}
\date{\today}
\address{Mathematisches Institut, Universit\"at zu K\"oln, Weyertal 86--90, D-50931 K\"oln, Germany}
\email{bruinier@math.uni-koeln.de }

\begin{abstract}
We give an introduction to the theory of Borcherds products and to some number theoretic and geometric applications. In particular, we discuss how the theory can be used to study the geometry of Hilbert modular surfaces.
\end{abstract}



\maketitle

\section{Introduction}
\label{sect1}

Infinite products play an important role in many branches of mathematics.
In number theory, they for instance provide an elegant way of encoding and manipulating combinatorial identities.
The product expansion of the generating function of the partition function is a well known example.
On the other hand, infinite products are a fundamental tool in complex analysis to construct meromorphic functions with prescribed zeros and poles,  the Weierstrass product theorem being a prominent example. 
In that way, they become interesting for the study of geometric problems.

In the first part of the present paper we will present some examples of particularly interesting infinite products, called Borcherds products, which are characterized by a striking symmetry property: They are modular forms for the orthogonal group of a suitable rational quadratic space of signature $(2,\ell)$. Although some very classical modular forms appear here, as for instance certain Eisenstein series or the $j$-function, most of these product expansions were only discovered rather recently by R.~Borcherds \cite{Bo1,Bo2}.

We will then consider more systematically the 
properties of Borcherds products on Hilbert modular surfaces. 
They can be used to recover important classical results on the geometry of such  surfaces. In addition, they provide a new approach to various arithmetic questions. 


Hilbert modular surfaces can be realized as arithmetic quotients associated to certain rational orthogonal groups of signature $(2,2)$. In particular, they are just very special instances of those Shimura varieties that can be obtained as quotients from orthogonal groups of signature $(2,\ell)$, the general setting of Borcherds' theory.
However, we feel that focusing on such a special case facilitates the presentation of central ideas. Moreover, the geometry of Hilbert modular surfaces is particularly beautiful.
Most of the results stated in sections \ref{sect:2}--\ref{sect:61} actually hold in greater generality.

The present paper is not intended to be a survey on Borcherds products. It only covers a small part of the many interesting aspects of the theory. For further expository articles, in particular for the connection to generalized Kac-Moody algebras,  we refer to  \cite{Bo6}, \cite{Bo5}, \cite{Bo4}, \cite{Ko}. 

I would like to thank E.~Freitag, W.~Kohnen, S.~Kudla, and U.~K\"uhn for many stimulating discussions and valuable suggestions.

\section{Infinite products and elliptic modular forms}
\label{sect:2}

Recall that an infinite product
\[
(1-a_1)(1-a_2)(1-a_3)\dots
\]
is said to converge absolutely, if the underlying series
\[
a_1+a_2+a_3+\dots
\]
converges absolutely. With this definition, an absolutely convergent infinite product vanishes, if and only if one
of its factors vanishes. So for instance  the product
\[
h(q)= \prod_{n=1}^\infty (1-q^n)
\]
converges absolutely for any complex number $q$ with $|q|<1$ and does not vanish there. Its underlying series is
the geometric series $\sum_{n\geq 1} q^n$.

This first example already has very interesting combinatorial and geometric properties. The function $1/h(q)$ is
holomorphic for $|q|<1$ and hence can be expanded in an infinite series
\begin{align*}
1/h(q) &= \sum_{n=0}^\infty p(n) q^n =1+q+2q^2+3 q^3+5 q^4+7q^5+11 q^6+15q^7+22 q^8+\dots.
\end{align*}
It is easily verified that the function $p(n)$ is the so called partition function. It counts the number of ways a
positive integer $n$ can be written as a sum of positive integers. The number of summands is unrestricted,
repetition is allowed, and the order of the summands is not taken into account. For instance there are $5$
partitions of $4$, namely $4$, $3+1$, $2+2$, $2+1+1$, $1+1+1+1$.

The arithmetic of partitions is quite involved and there are a number of unsolved questions about them (see e.g.~\cite{On}). One reason
is that the partition function grows rather rapidly. This is not apparent from the first few values give above,
but a quick computation (using e.g.~Maple) shows that for instance $p(100)=190\,569\,292$. In fact the celebrated
Hardy-Rademacher-Ramanujan asymptotic states
\begin{align}\label{partasym}
p(n)\sim \frac{e^{K\sqrt{n}}}{4n\sqrt{3}}
\end{align}
as $n\to\infty$, where $K=\pi\sqrt{2/3}$ (c.f.~\cite{Ap} chapter 5). This follows from the fact that $1/h(q)$
is closely related to a (weakly holomorphic) modular form of weight $-1/2$. The Fourier coefficients of such
modular forms satisfy similar asymptotics  in general.

In a different direction, if we put $q=e^{2\pi i\tau}$ for $\tau$ in the upper complex half plane $\H=\{\tau\in
\C;\; \Im(\tau)>0\}$, we are lead to the discriminant function
\begin{align}\label{defdelta}
\Delta(\tau)=q\cdot h(q)^{24}= q\prod_{n=1}^\infty (1-q^n)^{24}.
\end{align}
Since the product converges (locally uniformly) absolutely, it defines a nowhere vanishing holomorphic function on
$\H$. Moreover, $\Delta$ has a striking symmetry property with respect to the action of the modular group
$\Sl_2(\Z)$ on $\H$ by Moebius transformations $\tau\mapsto \kabcd \tau= \frac{a\tau+b}{c\tau+d}$. It satisfies
\[
\Delta\left(\frac{a\tau+b}{c\tau+d}\right)=(c\tau+d)^{12}\Delta(\tau),\qquad \abcd\in \Sl_2(\Z)
\]
(for a simple proof see \cite{Ko}).

To put this in a suitable context, we recall some basic notions of the theory of modular functions (see e.g.~\cite{Sh}, \cite{Ma}, \cite{Ap}
for more details). Let $k$ be an integer, $\Gamma$ a subgroup of finite index of  $\Gamma(1):=\Sl_2(\Z)$, and
$\chi$ an Abelian character of $\Gamma$. A function $f:\H\to \C$ is called a {\em weakly holomorphic  modular form} (of weight
$k$ with respect to $\Gamma$ and $\chi$), if
\begin{enumerate}
\item
$f\left(\frac{a\tau+b}{c\tau+d}\right)=\chi(\gamma)(c\tau+d)^{k}f(\tau)$ for all $\kabcd\in \Gamma$,
\item $f$ is holomorphic on $\H$,
\item $f$ is meromorphic at the cusps of $\Gamma$.
\end{enumerate}
We do not want to explain the last condition in too much detail. If $\Gamma=\Gamma(1)$ (which is all we need in
this section), then any function $f$ satisfying the first two conditions has a Fourier expansion
\[
f(\tau)= \sum_{n=-\infty}^{\infty} c(n) q^n,
\]
because $T=\kzxz{1}{1}{0}{1}\in \Gamma(1)$. Now the third condition means that only finitely many $c(n)$ with
$n<0$ are non-zero. If actually all $c(n)$ with $n<0$ (respectively $n\leq 0$) vanish, then $f$ is called a
{\em holomorphic modular form} (respectively {\em cusp form}) for $\Gamma(1)$. If $\Gamma$ has only finite index in
$\Gamma(1)$, one has to require similar conditions for finitely many $\Gamma(1)$-translates of $f$ corresponding
to the cusps of $\Gamma$, i.e., the $\Gamma$-orbits of $P^1(\Q)=\Q\cup\infty$.

We write $W_k$ (respectively $M_k$, $S_k$) for the space of weakly holomorphic modular forms (respectively
holomorphic modular forms, cusp forms) of weight $k$ for $\Gamma(1)$ with trivial character. One can show that
$M_k=\{0\}$ for $k<0$ or $k$ odd, $M_0=\C$, and $M_2=0$.

The properties of $\Delta$ stated above can be summarized by saying that $\Delta$ is a cusp form of weight $12$
for $\Gamma(1)$ with trivial character, i.e., an element of $S_{12}$.

Further examples of modular forms are provided by Eisenstein series. If $k$ is an even integer, $k>2$, we define
the Eisenstein series of weight $k$ for $\Gamma(1)$ by
\[
E_k(\tau)=\frac{1}{2}\sum_{\substack{c,d\in \Z\\\gcd(c,d)=1}} (c\tau+d)^{-k}.
\]
By comparison with a suitable integral, one checks that $E_k$ converges normally and defines a holomorphic function on
$\H$. Consequently the transformation behavior of a modular form of weight $k$ easily follows by reordering the
summation. The Fourier expansion of $E_k$ can be computed by means of the partial fraction expansion of the
cotangens. One finds that
\[
E_k(\tau)=1-\frac{2k}{B_k}\sum_{n\geq 1} \sigma_{k-1}(n)q^n,
\]
where $B_k$ denotes the $k$-th Bernoulli number and $\sigma_{k-1}(n)=\sum_{d\mid n}d^{k-1}$ the usual divisor sum
function. In particular $E_k$ is a (non-zero) holomorphic modular form of weight $k$ for $\Gamma(1)$. As a
consequence we get the decomposition $M_k=\C E_k\oplus S_k$ for $k>2$.

The non-vanishing of $\Delta$ on $\H$ implies that we have an isomorphism
\[
S_k\to M_{k-12},\quad f\mapsto f/\Delta.
\]
This observation is crucial. It can be deduced that the graded algebra $\bigoplus_{k\in \Z} M_k$ of modular forms
for $\Gamma(1)$ is a polynomial ring in $E_4$ and $E_6$.

To get a more geometric interpretation of modular forms for $\Gamma\subset \Gamma(1)$, one considers the modular
curve associated with $\Gamma$, that is, the quotient $Y_\Gamma=\Gamma\bs \H$. The complex structure of $\H$
induces a structure as a non-compact Riemann surface on $Y_\Gamma$. One obtains a compact Riemann surface
$X_\Gamma$ by adding finitely many points, namely the cusps of $\Gamma$. Every compact Riemann surface $X$ has a
natural structure as a projective algebraic variety over $\C$. The Zariski topology is given by taking for the
closed sets the finite subsets of $X$ together with $X$ itself. The structure sheaf $\calO$ is given by taking for
the regular functions $\calO(U)$ on an open subset $U\subset X$ all meromorphic functions on $X$, which are
holomorphic on $U$.

The modular curves $Y_\Gamma$ and $X_\Gamma$ associated to certain families of ``congruence subgroups'' $\Gamma$
can actually be defined over algebraic number fields and even over their rings of integers (see \cite{Sh}, \cite{DeRa}, \cite{DI}). This is due to their
moduli interpretation. For instance, the modular curve $\Gamma(1)\bs\H$ is the (coarse) moduli space of
isomorphism classes of elliptic curves over $\C$. In fact, by the theory of the Weierstrass $\wp$-function, an
elliptic curve over $\C$ is a torus, which in turn is given by a quotient $\C/L$, where $L\subset\C $ is a
lattice. Two elliptic curves $\C/L$ and $\C/L'$ are isomorphic, if and only if $L=aL'$ for some non-zero $a\in
\C-\{0\}$. In particular every isomorphism class has a representative $\C/L_\tau$, where $L_\tau=\Z+\tau\Z$ and $\tau
\in \H$. It is easily checked that two elliptic curves $\C/L_\tau$ and $\C/L_{\tau'}$ with $\tau,\tau' \in \H$ are
isomorphic, if and only if $\tau$ is equivalent to $\tau'$ with respect to the action of $\Gamma(1)$ on $\H$ (this
corresponds to a change of basis of the lattice). Hence the assignment $\tau\mapsto \C/L_\tau$ induces the
identification of $\Gamma(1)\bs\H$ with the moduli space. The point is that the moduli problem makes sense not
only over $\C$ but over any scheme over $\Z$, allowing to construct models over Dedekind rings.

We may regard $\Delta$ as a section of the line bundle of modular forms of weight $12$ on $\Gamma(1)\bs\H$. In
view of the moduli interpretation one checks that $\Delta$ assigns to the elliptic curve $\C/L_\tau$ its
discriminant $\Delta(\tau)$.

A further classical modular form is the $j$-function:
\[
j(\tau)=E_4^3(\tau)/\Delta(\tau)=q^{-1}+744+196884q+21493760q^2+\dots.
\]
Since $\Delta$ does not vanish on $\H$, the $j$-function is holomorphic on $\H$ . However, because of the presence
of the term $q^{-1}$ in the Fourier expansion, it is only meromorphic at the cusp $\infty$. Hence $j\in W_0$. In
terms of the moduli interpretation, $j$ assigns to the elliptic curve $\C/L_\tau$ its $j$-invariant $j(\tau)$. In
the theory of elliptic curves one shows that $j$ classifies elliptic curves over an algebraically closed field up
to isomorphism. Moreover, for every $z\in\C$ there is an elliptic curve with prescribed $j$-invariant $z$. In
other words, $j$ defines an analytic isomorphism
\[
\Gamma(1)\bs \H\longrightarrow \C.
\]
It extends to an isomorphism $X_{\Gamma(1)}\to P^1(\C)$ to the Riemann sphere $P^1(\C)$.

What does $j$ have to do with infinite products? There are different important infinite products involving $j$. Here
we only present one of these. To this end it is convenient to define $J(\tau)=j(\tau)-744\in W_0$. We denote the
Fourier expansion by $J(\tau)=\sum_{n\geq -1} c(n)q^n$ with coefficients $c(n)\in \Z$. In particular we have
$c(-1)=1$, and $c(0)=0$. We consider the modular form of two variables
\[
j(z_1)-j(z_2)
\]
of weight $0$ for the group $\Gamma(1)\times \Gamma(1)$. It vanishes at a point $(z_1,z_2)\in \H\times\H$, if and
only if $z_1=M z_2$ for some $M\in \Gamma(1)$, because of the injectivity of $j$. In other words, the divisor of
$j(z_1)-j(z_2)$ is given by quadratic equations
\begin{align}\label{hecke}
c z_1 z_2+d z_1-a z_2-b=0
\end{align}
with integral coefficients and $ad-bc=1$.

\begin{theorem}\label{thmjj}
The modular form $j(z_1)-j(z_2)$ has the infinite product expansion
\begin{align}\label{prodjj}
j(z_1)-j(z_2)=q_1^{-1}\prod_{\substack{m>0\\n\in \Z}}(1-q_1^m q_2^n)^{c(mn)}.
\end{align}
Here $q_j=e^{2\pi i z_j}$, and $c(n)$ is the $n$-th Fourier coefficient of $J(\tau)$. The product converges normally for
$\Im(z_1)\Im(z_2)>1$.
\end{theorem}

This beautiful identity was found in the 80's independently by Borcherds, Koike, Norton, and Zagier. It is the
denominator identity of the monster Lie algebra, which is a generalized Kac-Moody algebra with an action of the
Monster simple group. It is crucial in Borcherds' proof of the moonshine conjecture (see \cite{Bo4} for an
introduction and further references).

One may wonder why Theorem \ref{thmjj}  had not been discovered earlier. It only involves classical modular
functions that were well understood already at the end of the 19th century. One reason might be that the product
only converges on a sub-domain of $\H^2$. This is due to the fact that the $c(n)$, being the coefficients of a
{\em weakly} holomorphic modular form with a pole at the cusp $\infty$, grow rather rapidly:
\begin{align}\label{jasym}
c(n) \sim \frac{e^{4\pi\sqrt{n}}}{\sqrt{2} n^{3/4}},
\end{align}
as $n \to \infty$. This asymptotic is analogous to the one for the partition function \eqref{partasym}. The
restricted convergence of the product for $j(z_1)-j(z_2)$ implies that only part of the divisor can be read off
directly from the product. More precisely, it only tells us the vanishing along those divisors of type
\eqref{hecke} with $c=0$.

One might ask, whether the product expansions of the discriminant function $\Delta(\tau)$ and the function
$j(z_1)-j(z_2)$ have anything in common. This is in fact true. Both are particular examples of {\em Borcherds
products}. These are certain meromorphic modular forms in $\ell$ variables, which have a particular product
expansion, and arise as lifts of weakly holomorphic modular forms of weight $1-\ell/2$ for $\Gamma(1)$. Their
zeros and poles are explicitly given in terms of so-called  Heegner divisors (also referred to in the literature
as ``rational quadratic divisors'' or ``special divisors'').

In the next section we will make this informal definition more precise. Let us just remark here that
$j(z_1)-j(z_2)$ can be viewed as the Borcherds lift of the weight $0$ modular form $J(\tau)$, and  $\Delta(\tau)$
as the Borcherds lift of the classical weight $1/2$ Jacobi theta function $12\theta(\tau)=12\sum_{n\in \Z}
q^{n^2}$. It can be shown that the $j$-function and the Eisenstein series $E_k$, with $k=4,6,8,10,14$ are also
lifts of certain weakly holomorphic weight $1/2$ modular forms. For instance
\[
E_4=\prod_{n=1}^\infty(1-q^n)^{c(n^2)},
\]
where $c(n)$ denote the coefficients of the weight $1/2$ modular form (for $\Gamma_0(4)$ in the Kohnen plus-space)
\[
f=q^{-3}+4-240q+26760q^4-85995q^5+1707264 q^8-4096240q^9+44330496q^{12}+\dots,
\]
see \cite{Bo1}.

Finally, notice that Borcherds products naturally live on the hermitian symmetric space of type $IV$ associated to
the real orthogonal group $\Orth(2,\ell)$ of signature $(2,\ell)$. This hermitian space has complex dimension
$\ell$. They are constructed using a regularized theta lift for the dual reductive pair $(\Sl_2(\R),
\Orth(2,\ell))$. In small dimensions however, there are exceptional isomorphisms relating $\Orth(2,\ell)$ to other
classical Lie groups. For instance $\Orth(2,1)$ is essentially isomorphic to $\Sl_2(\R)$, and $\Orth(2,2)$ to
$\Sl_2(\R) \times \Sl_2(\R)$, which is implicitly used in the construction of the examples above. Moreover,
$\Orth(2,3)$ is  essentially isomorphic to the symplectic group $\Symp(2,\R)$ of genus $2$, and $\Orth(2,4)$ to
the hermitian symplectic group of genus $2$. In view of these isomorphisms, also the Heegner divisors alluded to
above become classically well known objects. For instance in the $\Orth(2,1)$ case, one gets Heegner points on
modular or Shimura curves, justifying the terminology. In the $\Orth(2,2)$ case, one can obtain Hirzebruch-Zagier
divisors on Hilbert modular surfaces, in the $\Orth(2,3)$ case Humbert surfaces on Siegel modular threefolds.

\section{Borcherds products on Hilbert modular surfaces}

We now want to generalize the $j(z_1)-j(z_2)$ example of the the previous section and study Borcherds'
construction of infinite automorphic products (in two variables) in a more systematic way. From a geometric point
of view, the underlying modular variety $\Gamma(1)^2\bs \H^2$ in that example was not very exciting. By means of
the $j$-function it is isomorphically mapped to the affine plane $\C^2$. If one wants to get more interesting
varieties one has to replace the discrete subgroup $\Gamma(1)\times\Gamma(1)\subset \Sl_2(\R)\times\Sl_2(\R)$ by
more complicated groups. An important family of discrete subgroups is provided by Hilbert modular groups of real
quadratic fields.

We first need to introduce some notation (see \cite{Fr}, \cite{Ge}, \cite{Go} for more details). Throughout we use
$z=(z_1,z_2)$ as a variable on $\H^2$ and write $(y_1,y_2)$ for its imaginary part. Let $K$ be the real quadratic
field of discriminant $D$. For simplicity we assume throughout that $D$ is a prime (hence $D\equiv 1\pmod{4}$ and
$K=\Q(\sqrt{D})$). We write $\OK$ for the ring of integers and $x\mapsto x'$ for the conjugation in $K$. The
Hilbert modular group $\GK=\Sl_2(\OK)$ associated with $K$ can be viewed as a discrete subgroup of $\Sl_2(\R)
\times\Sl_2(\R)$ by means of the two embeddings of $K$ into $\R$. In particular $\GK$ acts on $\H^2$ by
\[
(z_1,z_2)\mapsto \abcd(z_1,z_2)=\left( \frac{az_1+b}{cz_1+d},\, \frac{a'z_2+b'}{c'z_2+d'}\right).
\]
In the same way as with the modular curves $Y_\Gamma$ of the previous section, we consider the quotient
$Y_K=\GK\bs\H^2$, which has a structure as a non-compact complex surface. It can be compactified by adding $h_K$
points, namely the cusps of $\GK$, i.e., the $\GK$-orbits of $P^1(K)$. Here $h_K$ is the class number of $K$. In
contrast to the case of modular curves the resulting normal complex space $X_K$ is not regular. There are quotient singularities at the elliptic fixed points, and furthermore, the cusps are
highly singular points. 
By the theory of Hironaka the singularities can be resolved \cite{Hi}, there exists a desingularization $\pi: \widetilde{X}_K\to X_K$, such that $\D_K:=\pi^{-1}(X_K^{sing})$ is a divisor with normal crossings. 
The minimal resolution of singularities was constructed by Hirzebruch \cite{Hirzebruch}.

According to the theory of Baily-Borel, the complex spaces $Y_K$, $X_K$, and $\widetilde{X}_K$ all have natural
structures as quasi-projective algebraic varieties over $\C$. 
Moreover, $Y_K$ has a moduli interpretation as the (coarse) moduli space of isomorphism classes of Abelian surfaces with multiplication by
$\OK$ and a certain class of polarizations (see section \ref{sect:7}). 
This can be used to construct integral models.
The surfaces $Y_K$, $X_K$, and $\widetilde{X}_K$ are all  referred to
as {\em Hilbert modular surfaces} associated with  $K$.

On such surfaces there exist distinguished divisors, called Hirzebruch-Zagier divisors. (As already mentioned,
they play the role of the Heegner divisors in the informal definition of Borcherds products of the previous
section.) For every vector $(a,b,\lambda)$  of positive norm $ab-\lambda\lambda'$ in the lattice
$\Z^2\oplus\frakd^{-1}$, the subset
\[
(a,b,\lambda)^\perp=\{(z_1,z_2)\in \H^2; \quad az_1 z_2 +\lambda z_1 + \lambda' z_2 +b=0\}
\]
defines an analytic divisor on $\H^2$. Here $\frakd^{-1}=\frac{1}{\sqrt{D}}\OK$ is the inverse different of $K$. Observe the analogy with \eqref{hecke}.
Now let $m$ be a positive integer. The sum
\[
T(m)=\sum_{\substack{(a,b,\lambda)\in(\Z^2\oplus\frakd^{-1})/\{\pm 1\}\\
ab-\lambda\lambda'=m/D}} (a,b,\lambda)^\perp
\]
is a $\GK$-invariant analytic divisor on $\H^2$. It descends to an algebraic divisor on the quasi-projective
variety $Y_K$, the {\em Hirzebruch-Zagier divisor} of discriminant $m$. Moreover, we obtain Hirzebruch-Zagier divisors on $X_K$ by taking the closure of $T(m)$, and on $\widetilde{X}_K$ by taking the pullback with respect to the desingularization morphism.

One easily sees that $T(m)=\emptyset$, if and only if $\chi_D(m)=-1$, where $\chi_D$ is the quadratic character
corresponding to $K$ given by the Legendre symbol $\chi_D(x)=\left(\frac{D}{x}\right)$.  If $m$ is square-free and a norm of $\OK$, then the normalization of $T(m)$ on $X_K$ is isomorphic to the modular curve
$X_0(m)$.
If $m$ is square-free and not a norm of $\OK$ (but $\chi_D(m)=1$), then the normalization of $T(m)$ is isomorphic to a Shimura curve associated to a certain order in the
indefinite quaternion algebra $Q_m=\left(\frac{D,-m/D}{\Q}\right)$.

One can show that the Hirzebruch-Zagier divisors are actually defined over
$\Q$. Moreover, they often have a moduli interpretation, which can be used to extend them to integral models. For
instance, in the latter case, where the normalization of $T(m)$ is a Shimura curve, the points on $T(m)$ correspond to
Abelian surfaces with quaternionic multiplication by a certain order of $Q_m$.

Let $k$ be an integer. Recall that a meromorphic (respectively holomorphic) {\em Hilbert modular form} of weight
$k$ for the group $\Gamma_K$ is a meromorphic (respectively holomorphic) function $f$ on $\H^2$ satisfying the
transformation law
\begin{align}\label{eq:trafolaw}
f\left(\kabcd(z_1,z_2)\right)=(c z_1 +d)^{k}(c' z_2 +d')^{k} f(z_1,z_2)
\end{align}
for all $\kabcd\in \Gamma_K$. Notice that in contrast to the definition of modular forms for subgroups of
$\Sl_2(\Z)$ in the previous section, we do not have to require that $f$ be meromorphic (respectively holomorphic)
at the cusps. This is automatic by the Koecher principle.

\medskip

We are now ready to explain Borcherds' lift in the case of Hilbert modular surfaces. Let us first describe the
``input data'' which is used to feed the lift. We denote by $W_k(D,\chi_D)$ the space of weakly holomorphic
modular forms of weight $k$ for the Hecke group 
\begin{align}\label{gamma0}
\Gamma_0(D)=\left\{\abcd\in \Sl_2(\Z);\quad c\equiv 0\pmod{D}\right\}
\end{align} 
with character $\chi_D$
(where $\chi_D(\kabcd)=\chi_D(d)$).
Since $\kzxz{1}{1}{0}{1}\in \Gamma_0(D)$, 
any such modular form $f$ has a Fourier expansion of the form $f=\sum_{n\gg-\infty}c(n)q^n$. We let
$W_k^+(D,\chi_D)$ be the subspace of those $f\in W_k(D,\chi_D)$, whose Fourier coefficients $c(n)$ satisfy the so-called {\em plus space condition}, i.e., $c(n)=0$ whenever $\chi_D(n)=-1$. Moreover, we write $M^+_k(D,\chi_D)$ (respectively $S_k^+(D,\chi_D)$) for the subspace of
holomorphic modular forms (respectively cusp forms) in  $W_k^+(D,\chi_D)$.
For even $k\geq 2$, Hecke proved that there is a unique normalized 
Eisenstein series $E_k(\tau)$ in $M^+_k(D,\chi_D)$, and 
\begin{align}\label{eq:plussplit}
M^+_k(D,\chi_D)=S^+_k(D,\chi_D) \oplus \C E_k(\tau),
\end{align}
see \cite{He}, and section \ref{sect:61} here.

If $f=\sum_{n\in\Z} c(n) q^n\in \C((q))$ is a formal Laurent series, we put
\[
\tilde c(n)=\begin{cases}
c(n),& \text{if $n\not\equiv 0 \pmod{D}$,}\\
2c(n),&\text{if $n\equiv 0 \pmod{D}$.}
\end{cases}
\]
Modular forms in the plus space $W_k^+(D,\chi_D)$ can also be realized as vector valued modular forms for the
full modular group $\Gamma(1)$ transforming with a certain $D$-dimensional unitary representation (see \cite{BB}). As a
consequence, there is a bilinear pairing assigning to two modular forms $f\in W_k^+(D,\chi_D)$ and $g\in
W_{k'}^+(D,\chi_D)$ a modular form $\langle f, g\rangle\in W_{k+k'}$ for the  group $\Gamma(1)$. If $f=\sum_n
c(n)q^n$ and $g=\sum_n b(n)q^n$, it can be described in terms of the Fourier expansions as follows:
\begin{align}\label{pairfourier}
\langle f, g \rangle = \sum_{n\in\Z} \sum_{m\in \Z} \tilde c(m) b(Dn-m) q^n.
\end{align}
Notice that this pairing only depends on the transformation properties of modular forms in the plus space and
naturally extends to non-holomorphic forms.

\begin{theorem}[Borcherds]\label{hilbert}
Let $f=\sum_{n\gg -\infty}c(n)q^n$ be a weakly holomorphic modular form in $W_0^+(D,\chi_D)$ and assume that
$\tilde c(n)\in \Z$ for all $n<0$. Then there exists a meromorphic Hilbert modular form $\Psi(z_1,z_2,f)$ for $\GK$
(with some multiplier system of finite order) such that:
\begin{enumerate}
\item[(i)]
The weight of $\Psi$ is equal to the constant term $c(0)$ of $f$.
\item[(ii)]
The divisor $Z(f)$ of $\Psi$ is determined by the principal part of $f$ at the cusp $\infty$. It equals
\[Z(f)=\sum_{n<0} \tilde c(n)T(-n).\]
\item[(iii)]
Let $W\subset\H^2$ be a Weyl chamber attached to $f$ and put $N=\min \{n;\; c(n)\neq 0\}$. The function $\Psi$ has
the Borcherds product expansion
\[
\Psi(z_1,z_2,f)=q_1^{\rho} q_2^{\rho'} \prod_{\substack{\nu\in\frakd^{-1} \\ (\nu,W)>0}} \left(1-q_1^\nu q_2^{\nu'}
\right)^{\tilde c(D\nu\nu')},
\]
which converges normally for all $(z_1,z_2)$ with $y_1 y_2 > |N|/D$ outside the set of poles. Here $\rho\in K$ is
the Weyl vector corresponding to $W$ and $f$, and $q_j^\nu=e^{2\pi i \nu z_j}$ for $\nu\in K$.
%
%
\end{enumerate}
\end{theorem}

A few additional explanations regarding Weyl chambers are in order. For each $\lambda\in \frakd^{-1}$ of negative
norm the subset
\[
\lambda^\perp= \{(z_1,z_2)\in\H^2;\quad \lambda y_1 +\lambda' y_2=0\}
\]
is a hyperplane of real codimension $1$ in $\H^2$. Since $f$ is meromorphic at the cusps, and by reduction theory,
the union
\[
S(f)=\sum_{\substack{\lambda\in \frakd^{-1}\\ \lambda\lambda'<0}} c(D\lambda\lambda') \lambda^\perp
\]
has only finitely many $\GK$-orbits. It is invariant under the stabilizer in $\GK$ of the cusp $\infty$. Its
complement $\H^2-S(f)$ decomposes into connected components, which are called the Weyl chambers  attached to $f$.
To such a Weyl chamber $W$ (and $f$) one can associate the so-called Weyl vector $\rho\in K$, which we do not want
to define here (see \cite{BB} for its explicit computation in the present case, and \cite{Bo1, Bo2} for more
general facts). Moreover, if $\nu\in K$, then one writes $(\nu,W)>0$, if $\lambda y_1 +\lambda'y_2>0$ for all
$(z_1,z_2)\in W$.

If $\dv(\Psi)$ is compact in $Y_K$, i.e., if the normalizations of the irreducible components are Shimura curves, then
$S(f)$ is empty and there is just the one Weyl chamber $W=\H^2$. In this case $\rho=0$, and the condition
$(\nu,W)>0$ becomes just the condition that $\nu$ be totally positive.

Theorem \ref{hilbert} is contained (in a slightly different formulation) in Theorem 13.3 of \cite{Bo2} (to obtain
the above form see \cite{BB}). The idea of the proof is as follows.

First, we notice that by an elementary argument the Fourier coefficients of $f$ are bounded by
\begin{align}\label{est}
c(n)=O\left(e^{4\pi \sqrt{|Nn|}}\right),\qquad n\to \infty,
\end{align}
see \cite{BF} section 3. This implies the convergence of the product in the stated region. The estimate
\eqref{est} is also a consequence of the (much more precise) Hardy-Rademacher-Ramanujan asymptotic for the
coefficients of weakly holomorphic modular forms, the general theorem behind \eqref{partasym} and \eqref{jasym}.

As already mentioned, the group $(\Sl_2(\R)\times \Sl_2(\R))/\{\pm 1\}$ is isomorphic to the connected component
of real orthogonal group $\Orth(2,2)$. The latter group and $\Sl_2(\R)$ form a dual reductive pair in the sense of
Howe \cite{Ho}. Thus we can construct Hilbert modular forms for $\GK$ from modular forms on $\Sl_2(\R)$ by
integrating against a certain kernel function, the Siegel theta function $\Theta_S(\tau;z_1,z_2)$ associated to
the lattice $L=\Z^2\oplus\OK$ of signature $(2,2)$. Here the Siegel theta function is a certain non-holomorphic
modular form for $\Gamma_0(D)$ satisfying the plus-space condition in the variable $\tau$, and a
$\GK$-invariant function in the variable $(z_1,z_2)$. This suggests that we look at the integral
\begin{align}\label{thetaint}
\int_\calF \left\langle f(\tau),\Theta_S(\tau;z_1,z_2)\right\rangle\,d\mu.
\end{align}
Here $\calF=\{\tau\in \H;\; |\tau|\geq 1, \;|u|\leq 1/2\}$ denotes the standard fundamental domain for the action of
$\Gamma(1)$ on $\H$, and $d\mu =\frac{du\,dv}{v^2}$  is the invariant measure on $\H$ (with $\tau=u+iv$).
Since $f$ and the Siegel theta function are in weight $0$, the integrand is  $\Gamma(1)$-invariant, so that the
integral makes formally sense. Unfortunately, because $f$ grows exponentially as $\tau$ approaches the cusp
$\infty$, it diverges wildly. However, Harvey and Moore discovered that it is possible to regularize the integral
\cite{HM}. Following their idea, Borcherds defines the regularized theta lift $\Phi(z_1,z_2,f)$ of $f$ to be the
constant term in the Laurent expansion at $s=0$ of the meromorphic continuation in $s$ of
\begin{align}\label{regint}
\lim_{t\to \infty} \int_{\calF_t} \left\langle f(\tau),\Theta_S(\tau;z_1,z_2)\right\rangle v^{-s}\,d\mu.
\end{align}
Here $\calF_t=\{\tau\in \calF; \;|v|\leq t\}$ is the truncated fundamental domain. One can show that this regularized
integral still makes sense even though \eqref{thetaint} does not. It defines a $\Gamma_K$-invariant real analytic
function on $\H^2- \supp(Z(f))$ with a logarithmic singularity\footnote{If $X$ is a normal complex space, $D\subset X$ a Cartier divisor, and $f$ a smooth function on $X-\supp(D)$, then $f$ has a logarithmic singularity along $D$, if for any local equation $g$ for $D$ on an open subset $U\subset X$, the function $f-\log|g|$ is smooth on $U$.} 
along the divisor $-4 Z(f)$.

Moreover, it can be shown that the limit in \eqref{regint} exists and is holomorphic at $s=0$, if the constant
term $c(0)$ of $f$ vanishes. It follows that $\Phi(z_1,z_2,f)$ is equal to
\[
\lim_{t\to \infty} \int_{\calF_t} \big(\left\langle f(\tau),\Theta_S(\tau;z_1,z_2)\right\rangle -\tilde c(0)
v\big)\, d\mu + A\tilde c(0),
\]
where $A$ is the constant term in the Laurent expansion at $s=0$ of $\lim_{t\to\infty} \int_{\calF_t} v^{1-s}
d\mu$. This could be taken as an alternative definition.

The Fourier expansion of $\Phi(z_1,z_2,f)$ can be computed explicitly by applying some partial Poisson summation
on the theta kernel. It turns out that
\[
\Phi(z_1,z_2,f)= -2\log\left| \Psi(z_1,z_2,f)^2 (16\pi^2 y_1 y_2)^{c(0)} \right| + 2c(0)\left(\log(8\pi)-\Gamma'(1)\right),
\]
giving the meromorphic continuation, the divisor, and the transformation behavior of the infinite product
$\Psi(z_1,z_2,f)$.

\medskip

Notice that a weakly holomorphic modular form $f=\sum_{n} c(n) q^n$  in $W_k^+(D,\chi_D)$ of weight $k\leq 0$ is uniquely determined by
its {\em principal part}
\[
\sum_{n<0} c(n)q^n\in \C[q^{-1}].
\]
For the difference of two elements of $W_k^+(D,\chi_D)$ with the same principal part is a holomorphic modular form
of weight $k\leq 0$ with Nebentypus, hence vanishes identically. Moreover, the Galois group $\Gal(\C/\Q)$ acts on weakly
holomorphic modular forms by mapping $f$ to
\[
f=\sum_{n} c(n)^\sigma  q^n,\qquad \sigma\in \Gal(\C/\Q).
\]
Here $ c(n)^\sigma $ denotes the Galois conjugate of $c(n)$. This follows from the fact that the spaces
$W_k(D,\chi_D)$ have bases of modular forms with integral rational Fourier coefficients (see \cite{DeRa}, \cite{DI}).
Consequently, if $f$ has principal part in $\Q[q^{-1}]$, then all its Fourier coefficients are rational with
bounded denominators.
We may conclude that some power of any holomorphic Borcherds product has coprime integral rational Fourier coefficients and
trivial multiplier system. This observation is crucial. By the $q$-expansion principle (see \cite{Ge}, chapter X.3
Theorem 3.3), such a modular form corresponds to a section of the line bundle of Hilbert  modular forms over $\Z$
on the moduli stack over $\Z$ representing the moduli problem ``isomorphism classes of Abelian schemes with real
multiplication by $\OK$''. Hence Borcherds products can be viewed as modular forms over $\Z$. In fact, they
provide a powerful tool to study arithmetic intersection numbers on Hilbert modular surfaces, see \cite{BBK} and section \ref{sect:8} here.

\medskip

We conclude this section with an example for Theorem \ref{hilbert}. We consider the Hilbert modular group $\GK$ of
the real quadratic field $K=\Q(\sqrt{5})$. The fundamental unit of $K$ is $\eps_0=\frac{1}{2}(1+\sqrt{5})\in \OK$.
Gundlach constructed a particular holomorphic Hilbert modular form $\Theta(z_1,z_2)$ of weight $5$ for $\GK$ as a
product of $10$ theta functions  of weight $1/2$ \cite{Gu}. He showed that the divisor of $\Theta$ is equal to
$T(1)$ and used this fact to determine the graded algebra of Hilbert modular forms for $\GK$. From the
construction one also finds that $\Theta$ has Fourier coefficients in $\Z$ with greatest common divisor $64$.

One can recover Gundlach's function using the Borcherds lift: We need to look at the ``input space''
$W_0^+(5,\chi_5)$. Using some basic facts on modular forms for $\Gamma_0(5)$ due to Hecke one finds that there is
a modular form $f_1\in W_0^+(5,\chi_5)$ with Fourier expansion
\[
f_1=q^{-1} + 5 + 11\,q - 54\,q^{4} + 55\,q^{5} + 44\,q^{6} - 395\,q^{9} + 340\,q^{10} + 296\,q^{11} -
1836\,q^{14}+\dots.
\]
If we plug this into the Borcherds lift, we get a Hilbert modular form $\Psi$ for $\GK$ of weight $5$ with divisor
$T(1)$. Hence $\Psi$ must be a constant multiple of $\Theta$. From the Borcherds product expansion it follows that
$\Psi$ has coprime Fourier coefficients in $\Z$. Consequently the constant factor is $1/64$. If we compute the
Weyl vector as in \cite{BB}, we obtain the product expansion
 \begin{equation}\label{gund}
\frac{1}{64}\Theta(z_1,z_2)=q_1^{\eps_0/\sqrt{5}}q_2^{(\eps_0/\sqrt{5})'} \prod_{\substack{\nu\in\frakd^{-1} \\
\eps_0\nu'-\eps_0'\nu>0}} \left(1-q_1^\nu q_1^{\nu'}\right)^{\tilde c(5\nu\nu')},
\end{equation}
where the $c(n)$ denote the Fourier coefficients of $f_1$.

\section{Obstructions for constructing Borcherds products}

Borcherds' theorem of the previous section provides a way of constructing many Hilbert modular forms with known
divisor supported on Hirzebruch-Zagier divisors $T(m)$. It is natural to seek for a precise description of those
linear combinations of Hirzebruch-Zagier divisors, which are divisors of Borcherds products. Since the divisor of
a Borcherds product $\Psi(z_1,z_2,f)$ is determined by the principal part of the weakly holomorphic modular form
$f$, which is used to construct it, it suffices to understand which Fourier polynomials $\sum_{n<0} c(n)q^{n}\in
\C[q^{-1}]$ can occur as principal parts of elements of $W_0^+(D,\chi_D)$. A necessary condition is easily
obtained. If $f\in W_k^+(D,\chi_D)$ with Fourier coefficients $c(n)$, and $g\in M_{2-k}^+(D,\chi_D)$ with Fourier
coefficients $b(n)$, then the pairing $\langle f,g\rangle$ is a weakly holomorphic modular form of weight $2$ for
$\Gamma(1)$. Thus
\[
\langle f,g\rangle d\tau
\]
is a meromorphic differential on the Riemann sphere whose only pole is at the cusp $\infty$. By the residue
theorem its residue has to vanish. But this residue is just the constant term in the Fourier expansion of $\langle
f,g\rangle$. In view of \eqref{pairfourier} we find that
\begin{align}\label{neccond}
\sum_{n\leq 0} \tilde c(n) b(-n)=0.
\end{align}

Using Serre duality for vector bundles on Riemann surfaces, Borcherds showed that this condition is essentially also
sufficient (see \cite{Bo3} and \cite{BB} Theorem 6).

\begin{theorem}\label{serre}
There exists a weakly holomorphic modular form $f\in W_k^+(D,\chi_D)$ with prescribed principal part
$\sum_{n<0}c(n)q^n$ (where $c(n)=0$ if $\chi_D(n)=-1$), if and only if
\[
\sum_{n<0} \tilde c(n)b(-n)=0
\]
for every cusp form $g=\sum_{m>0} b(m)q^m$ in $S_{2-k}^+(D,\chi_D)$.
\end{theorem}

This result shows that $S_2^+(D,\chi_D)$ is precisely the space of obstructions for constructing Borcherds
products on $Y_K$ with prescribed divisor. (In the same way $M_2^+(D,\chi_D)$ can be viewed as the space of obstructions for constructing Borcherds
products with prescribed divisor and weight.) 
The dimension of this space can be computed by means of the Riemann-Roch theorem or the Selberg
trace formula. In our case, where $D$ is prime, this was already done by Hecke \cite{He}. One finds that
\[
\dim S_2^+(D,\chi_D) = \dim M_2^+(D,\chi_D)-1= \left[\frac{D-5}{24}\right].
\]
In particular $S_2^+(D,\chi_D)=\{0\}$ for the primes $D=5, 13, 17$. In these cases there are no obstructions, and
for any $T(m)$ there is a Borcherds products product with divisor $T(m)$. For all other prime discriminants $D$
there are obstructions. Then for instance $T(1)$ is not the divisor of a Borcherds product, since there is a
normalized Hecke eigenform $g$ in $S_2(D,\chi_D)$. Its first Fourier coefficient is $1$, and the ``real part'' of $g$
is an element of $S_2^+(D,\chi_D)$, whose first Fourier coefficient is equal to $1$ as well (see \cite{Ge},
chapter VI.4). Hence $q^{-1}$ cannot be the principal part of a weakly holomorphic modular form in
$W_0^+(D,\chi_D)$.

Let us look at Theorem \ref{serre} from a little different angle. Let $\C[q^{-1}]^+$ (respectively $\C[[q]]^+$) be
the space of Fourier polynomials (respectively formal power series), whose coefficients satisfy the plus space
condition. We define a non-degenerate bilinear pairing between 
these spaces by putting
\[
\{f,g\}=\sum_{n\leq 0} \tilde c(n)b(-n)
\]
for $f=\sum_{n\leq 0}c(n)q^n\in\C[q^{-1}]^+$ and $g=\sum_{m\geq 0}b(m)q^m\in\C[[q]]^+$. 
For $k\leq 0$ the space $W_k^+(D,\chi_D)$ is
identified with a subspace of $\C[q^{-1}]^+$ by mapping a weakly holomorphic modular form with Fourier coefficients $c(n)$ to $\sum_{n\leq 0}c(n)q^n$.
The space $M_{2-k}^+(D,\chi_D)$ is identified with a subspace of $\C[[q]]^+$ by mapping a modular form to its
$q$-expansion. Now Theorem \ref{serre} implies that the orthogonal complement of $S_{2-k}^+(D,\chi_D)$ with respect to the pairing $\{\cdot,\cdot\}$ is equal to $W_k^+(D,\chi_D)\oplus \C$. Using the splitting \eqref{eq:plussplit} of $M_{2-k}^+(D,\chi_D)$  one concludes that the orthogonal complement of $M_{2-k}^+(D,\chi_D)$ is precisely $W_k^+(D,\chi_D)$. Since the pairing is non-degenerate, and since $M_{2-k}^+(D,\chi_D)$ has finite dimension, it follows by linear algebra that conversely
$M_{2-k}^+(D,\chi_D)$ is the orthogonal complement of $W_k^+(D,\chi_D)$. In other words:

\begin{corollary}\label{serre2}
A formal power series $\sum_{m\geq 0} b(m)q^m\in \C[[q]]^+$ is the $q$-expansion of a modular form in
$M_{2-k}^+(D,\chi_D)$, if and only if
\[
\sum_{n\leq 0} \tilde c(n)b(-n)=0
\]
for every  $f=\sum_{n} c(n)q^n$ in $W_{k}^+(D,\chi_D)$.
\end{corollary}

Since the pairing $\{\cdot,\cdot\}$ is defined over $\Q$ with respect to the natural rational structures
$\Q[q^{-1}]^+$ and $\Q[[q]]^+$, and since $M_{2-k}^+(D,\chi_D)$ and $W_{k}^+(D,\chi_D)$ have bases of modular
forms with integral coefficients, an analogous assertion holds for modular forms over $\Q$. Moreover, it suffices to
check the condition in this corollary for every $f\in W_{k}^+(D,\chi_D)$ with integral coefficients.

\medskip

If $X$ is a regular projective algebraic variety, we write $\ch^1(X)$ for its first Chow group, i.e., the
group of algebraic divisors on $X$ modulo rational equivalence. Furthermore, we put $\ch^1(X)_\Q=\ch^1(X)\otimes_\Z \Q$. 
Recall that $\ch^1(X)$ is isomorphic to the Picard group of $X$, the group of isomorphism classes of algebraic line bundles on $X$. The isomorphism is given by mapping a line bundle $\calL$ to the class $\cc_1(\calL)$ of the divisor of a rational section of $\calL$.
The Chow group $\ch^1(X)$ is an important invariant of $X$. It is finitely generated.

Meromorphic (respectively holomorphic) Hilbert modular forms can be interpreted as rational (respectively regular) sections of the sheaf $\calM_k(\C)$ of modular forms, which can be defined as follows: If we write $p:\H^2\to Y_K$ for the canonical projection, then the sections over an open subset $U\subset\Gamma\bs \H^2$ are holomorphic functions on $p^{-1}(U)$, satisfying the transformation law \eqref{eq:trafolaw}. By the Koecher principle, this sheaf on $Y_K$ extends to $X_K$. Moreover, we obtain a sheaf on $\widetilde{X}_K$, also denoted by $\calM_k(\C)$, by taking the pullback with respect to the desingularization morphism.
By the theory of Baily-Borel, there is a positive integer $n_K$ such that $\calM_k(\C)$ is an algebraic line bundle if $n_K | k$, and therefore defines an element of $\Pic(\widetilde{X}_K)$. 
Notice that $\calM_{nk}(\C)=\calM_{k}(\C)^{\otimes n}$ for any positive integer $n$.
If $k$ is any rational number, we chose an integer $n$ such that $nk$ is a positive integer divisible by $n_K$ and put $\cc_1(\calM_k(\C))=\frac{1}{n} \cc_1(\calM_{nk}(\C))\in \ch^1(\widetilde{X}_K)_\Q$.

It is natural to study the positions of the
Hirzebruch-Zagier divisors in $\ch^1(\widetilde{X}_K)$. 
To this end we consider the generating series
\begin{align}\label{genseries}
A(\tau) = 
\cc_1(\calM_{-1/2}(\C)) + \sum_{m>0} T(m) q^m \in \Q[[q]]^+\otimes_\Q \ch^1(\widetilde{X}_K)_\Q.
\end{align}
Combining Theorem \ref{hilbert} and Corollary \ref{serre2} one obtains the following striking application.

\begin{theorem}\label{hirzebruchzagier}
The Hirzebruch-Zagier divisors generate a subspace of $\ch^1(\widetilde{X}_K)_\Q$ of dimension $\leq \dim(M_{2}^+(D,\chi_D))$.
The generating series $A(\tau)$ is a modular form in $M_{2}^+(D,\chi_D)$ with values in
$\ch^1(\widetilde{X}_K)_\Q$, i.e., an element
of
$M^+_2(D,\chi_D)\otimes_\Q\ch^1(\widetilde{X}_K)_\Q$. 
\end{theorem}

In other words, if $\lambda$ is a linear functional on $\ch^1(\widetilde{X}_K)_\Q$, then
\begin{align*}
 \lambda\left(\cc_1(\calM_{-1/2}(\C))\right) + \sum_{m>0}  \lambda(T(m)) q^m \in M_{2}^+(D,\chi_D).
\end{align*}
A typical linear functional, one can take for $\lambda$, is given by the intersection pairing 
with a divisor on $\widetilde{X}_K$. Theorem \ref{hirzebruchzagier} was first proved by Hirzebruch and Zagier 
by
computing intersection numbers on $\widetilde{X}_K$ of Hirzebruch-Zagier divisors with other such divisors and
with the exceptional divisors coming from the resolution of the cusp singularities \cite{HZ}. Their discovery triggered
important investigations by several people, showing that more generally periods of certain special cycles in
arithmetic quotients of orthogonal or unitary type can be viewed as the coefficients of Siegel modular forms. For
instance, Oda considered cycles on quotients of $\Orth(2,\ell)$ given by embedded quotients of $\Orth(1,\ell)$
\cite{Od}, and Kudla-Millson studied more general cycles on quotients of $\Orth(p,q)$ and $\Uni(p,q)$ using the
Weil representation and theta functions with values in closed differential forms \cite{KM1,KM2,KM3}, see also \cite{Fu} for the case of non-compact quotients. The
relationship of the Kudla-Millson lift and Borcherds' regularized theta lift is clarified in \cite{BF}.

Using Borcherds products, Theorem \ref{hirzebruchzagier} can be proved as follows (see \cite{Bo3}). In view of
Corollary \ref{serre2} it suffices to show that
\[
\tilde c(0)\cc_1(\calM_{-1/2}(\C)) + \sum_{n<0} \tilde c(n) T(-n)=0 \in \ch^1(\widetilde{X}_K)_\Q
\]
for every  $f=\sum_{n} c(n)q^n$ in $W_{0}^+(D,\chi_D)$ with integral Fourier coefficients. But this is an immediate
consequence of Theorem \ref{hilbert}: The Borcherds lift of $f$ is a rational section of $\calM_{c(0)}(\C)$ with divisor $\sum_{n<0} \tilde c(n) T(-n)$.
Notice that we have only used (i) and (ii) of Theorem \ref{hilbert}. The product expansion (iii) is not required.
Finally, we mention that this argument generalizes to Heegner divisors on
quotients of $\Orth(2,\ell)$.

\section{Converse theorems}

By Theorem \ref{serre} of the previous section we know precisely which linear combinations of Hirzebruch-Zagier
divisors occur as divisors of Borcherds products on $Y_K$. Here it is natural to ask, whether every Hilbert modular
form on $Y_K$, whose divisor is a linear combination of Hirzebruch-Zagier divisors, is a Borcherds product, i.e.,
in the image of the lift of Theorem \ref{hilbert}. In this section we discuss this question in some detail. To
answer it, we first simplify the problem. We extend the Borcherds lift to a larger space of ``input modular
forms'', given by certain Maass wave forms, and answer the question for this extended lift. In that way we are led
to automorphic Green functions associated with Hirzebruch-Zagier divisors.

Let $k$ be an integer, $\Gamma$ a subgroup of finite index of  $\Gamma(1)$, and $\chi$ a character of $\Gamma$. A
twice continuously differentiable function $f:\H\to \C$ is called a {\em weak Maass form} (of weight $k$ with
respect to $\Gamma$ and $\chi$), if
\begin{enumerate}
\item
$f\left(\frac{a\tau+b}{c\tau+d}\right)=\chi(\gamma)(c\tau+d)^{k}f(\tau)$ for all $\kabcd\in \Gamma$,
\item $f$ has polynomial growth at the cusps of $\Gamma$ (in terms of local parameters),
\item $\Delta_k f(\tau)=0$.
\end{enumerate}
Here
\begin{equation}\label{deflap}
\Delta_k = -v^2\left( \frac{\partial^2}{\partial u^2}+ \frac{\partial^2}{\partial v^2}\right) + ikv\left(
\frac{\partial}{\partial u}+i \frac{\partial}{\partial v}\right)
\end{equation}
denotes the usual hyperbolic Laplace operator in weight $k$ and $\tau=u+iv$.

So if we compare this with the definition of a weakly holomorphic modular form, we see that we simply
replaced the condition that $f$ be holomorphic on $\H$ by the weaker condition that $f$ be annihilated by
$\Delta_k$, and the meromorphicity at the cusps by the corresponding growth condition. The third condition implies
that $f$ is actually real analytic. Because of the transformation behavior, it has a Fourier expansion involving
besides the exponential function a second type of Whittaker function. (See \cite{BF} section 3 for more details.)

There are two fundamental differential operators on modular forms for $\Gamma$, the Maass raising and lowering
operators \begin{align*} R_k  =2i\frac{\partial}{\partial\tau} + k v^{-1} \qquad \text{and} \qquad L_k  = -2i v^2
\frac{\partial}{\partial\bar{\tau}}.
\end{align*}
If $f$ is a differentiable function on $\H$ satisfying the transformation law (1) above in weight $k$, then $L_k
f$ transforms in weight $k-2$, and $R_k f$ in weight $k+2$. It can be shown that the assignment
\[
f(\tau)\mapsto \xi_k(f)(\tau):=v^{k-2} \overline{L_k f(\tau)} = R_{-k} v^k\overline{ f(\tau)}
\]
defines an antilinear map $\xi_k$ from weak Maass forms of weight $k$ to weakly holomorphic modular forms of
weight $2-k$. Its kernel is precisely the space of weakly holomorphic modular forms in weight $k$.

We write $N_k(D,\chi_D)$ for the space of weak Maass forms of weight $k$ with respect to $\Gamma_0(D)$ and
$\chi_D$. Let us have a closer look at map $\xi_k: N_k(D,\chi_D) \to W_{2-k}(D,\chi_D)$. We denote by
$\calN_k(D,\chi_D)$ the inverse image of $S_{2-k}(D,\chi_D)$ under $\xi_k$, and its plus subspace by
$\calN_k^+(D,\chi_D)$. (Note that our notation is not consistent with the notation of \cite{BF}.)

\begin{theorem}\label{exseq}
We have the following exact sequence:
\[
\xymatrix{ 0\ar[r]& W_k^+(D,\chi_D) \ar[r]& \calN_k^+(D,\chi_D) \ar[r]^{\xi_k}& S_{2-k}^+(D,\chi_D) \ar[r] & 0. }
\]
\end{theorem}

This can be proved using Serre duality for the Dolbeault resolution of the structure sheaf on a modular curve (see
\cite{BF} Theorem 3.7) or by means of Hejhal-Poincar\'e series (see \cite{Br2} chapter 1).

For every weak Maass form $f\in \calN_k^+(D,\chi_D)$ there is a unique Fourier polynomial $P(f)=\sum_{n<0}c(n)q^n$
in $\C[q^{-1}]^+$ such that $f(\tau)-P(f)(\tau)$ is bounded as $v\to\infty$, uniformly in $u$. It is called the
principal part of $f$. This generalizes the notion of the principal part of a weakly holomorphic modular form. One
can show that every prescribed Fourier polynomial as above occurs as the principal part of a unique $f\in
\calN_k^+(D,\chi_D)$. This is a key fact, which suggests to study the Borcherds lift of weak Maass forms.

If $f\in \calN_0^+(D,\chi_D)$, then we define its regularized theta lift $\Phi(z_1,z_2,f)$ by \eqref{regint}, in
the same way as for weakly holomorphic modular forms. One can show that the regularized theta integral defines a
$\Gamma_K$-invariant function on $\H^2$ with a logarithmic singularity along $-4Z(f)$, where
\[
Z(f)=\sum_{n<0} \tilde c(n)T(-n),
\]
and $\sum_{n<0}c(n)q^n$ denotes the principal part of $f$ \cite{Br2}, \cite{BF}. It is almost harmonic (outside the singularities) in
the following sense. If $\Delta^{(1)}$ and $\Delta^{(2)}$ denote the $\Sl_2(\R)$-invariant hyperbolic Laplace
operators on $\H^2$ in the first and second variable, then $\Delta^{(j)}\Phi(z_1,z_2,f)=\text{constant}$.

The Fourier expansion of $\Phi(z_1,z_2,f)$ can be computed explicitly. It can be used to determine the growth
behavior at the boundary of $Y_K$ in $\widetilde{X}_K$. It turns out that the boundary singularities are of $\log$
and $\log$-$\log$ type. More precisely, one can view $\frac{1}{4}\Phi(z_1,z_2,f)$ as 
a pre-log-log Green function for the divisor
$Z(f)$ on $\widetilde{X}_K$
in the sense of \cite{BKK} (see section \ref{sect:7} here, and \cite{BBK} Proposition 2.15).

Moreover, one finds that $\Phi(z_1,z_2,f)$ can be split into a sum
\begin{align}\label{split}
\Phi(z_1,z_2,f)=-2 \log| \Psi(z_1,z_2,f)|^2 + \xi(z_1,z_2,f),
\end{align}
where $\xi(z_1,z_2,f)$  is real analytic on the whole domain $\H^2$ and $\Psi(z_1,z_2,f)$ is a meromorphic
function on $\H^2$ whose divisor equals $Z(f)$. If $f$ is weakly holomorphic, the function $\xi(z_1,z_2,f)$ is 
simply equal to $2c(0)\left(\log(8\pi)-\Gamma'(1)- \log(16\pi^2 y_1 y_2)\right)$ and we are back in the case of Borcherds' original lift. However, if $f$ is
an honest weak Maass form, then $\xi$ is a complicated function and $\Psi$ far from being modular.

In any case, via the usual Poincar\'e-Lelong argument, the above splitting implies that the $(1,1)$ form
\begin{equation}\label{deflb}
\Lambda_B(f):=\frac{1}{4} dd^c \Phi (z_1,z_2,f)
= \frac{1}{4} dd^c \xi(z_1,z_2,f)
\end{equation}
represents the Chern class of the divisor $Z(f)$ in the second cohomology $H^2(Y_K)$. One can further show that it
is a square integrable harmonic representative. (It can also be regarded as a pre-log-log form on $\widetilde{X}_K$, representing the class of $T(m)$ on $\widetilde{X}_K$ in $H^2(\widetilde{X}_K)$.)

Using suitable $L^p$-estimates for the functions $\Phi(z_1,z_2,f)$, and results of Andreotti-Vesentini and Yau on
\mbox{(sub-)} harmonic functions on complete Riemann manifolds that satisfy such integrability conditions, the
following {\em weak converse theorem} for the Borcherds lift can be proved (see \cite{Br2} chapter 5).

\begin{theorem}\label{converse1}
Let $F$ be a meromorphic Hilbert modular form of weight $r$ for $\Gamma_K$, whose divisor $\dv(F)=\sum_{n<0}
\tilde c(n)T(-n)$ is a linear combination of Hirzebruch-Zagier divisors. Then
\[
-2 \log | F(z_1,z_2)^2 (16\pi^2 y_1y_2)^r| = \Phi(z_1,z_2,f)+\text{constant},
\]
where $f$ is the unique weak Maass form in $\calN_0^+(D,\chi_D)$ with principal part $\sum_{n<0}c(n)q^n$.
\end{theorem}

\begin{corollary}
The assignment $Z(f)\mapsto \Lambda_B(f)$ defines a linear map from the subgroup $\ch^1_{HZ}(\widetilde{X}_K)$ of
$\ch^1(\widetilde{X}_K)$, generated by the Hirzebruch-Zagier divisors, to $\calH^{1,1}(Y_K)$, the space of square integrable
harmonic $(1,1)$-forms on $Y_K$.
\end{corollary}

Summing up, we get the following commutative diagram:
\begin{align}\label{eq:diag}
\xymatrix{
\calN_0^+(D,\chi_D) \ar[r] \ar[d]
& \calN_0^+(D,\chi_D)/W_0^+(D,\chi_D)\ar[d]\ar[r]^-{\sim}_-{\xi_0}& S_2^+(D,\chi_D)\\
\z^1_{HZ}(\widetilde{X}_K) \otimes_\Z \C \ar[r] & \ch^1_{HZ}(\widetilde{X}_K) \otimes_\Z\C \ar[r]
&  \calH^{1,1}(Y_K) }.
\end{align}
Here $\z^1_{HZ}(\widetilde{X}_K)$ denotes the subgroup of the divisor group $\z^1(\widetilde{X}_K)$ generated by the $T(m)$ ($m\in \N$).
The left vertical arrow is defined by $f\mapsto Z(f)-\tfrac{c(0)}{k_G}\dv{G}$, where $c(0)$ denotes the constant term of $f$, and $G$ a fixed meromorphic Borcherds product of weight $k_G$. The vertical arrow in the middle is given by $f\mapsto Z(f)-c(0)\cc_1(\calM_1(\C))$.

In particular, the above diagram gives rise to a linear map $S_2^+(D,\chi_D)\to
\calH^{1,1}(Y_K) $. It can be explicitly described in terms of the Fourier expansions. 
One finds that
the image is in the subspace of $\calH^{1,1}(Y_K)$ given by forms which are symmetric with respect to
the interchange of the coordinates $z_1,z_2$. It is known that this subspace is isomorphic to $\C\Omega\oplus S_2(\GK)$, where 
\begin{align}\label{defomega}
\Omega=\frac{dx_1\,dy_1}{4\pi y_1^2}+\frac{dx_2\,dy_2}{4\pi y_2^2}
\end{align}
is the symmetric invariant K\"ahler form on $\H^2$, and $S_2(\GK)$ denotes the space
of Hilbert cusp forms of weight $2$ for $\GK$.
Consequently, we get a linear map $S_2^+(D,\chi_D)\to S_2(\GK)$.

To answer the surjectivity question for the Borcherds lift raised at the beginning of this section, it now
suffices to show that this map is injective. This can for instance be deduced by means of the Fourier expansion of
the image. We obtain the following  {\em strong converse theorem} for the Borcherds lift (see \cite{Br1},
\cite{Br2} chapter 5).

\begin{theorem}\label{converse2}
Let $F$ be a meromorphic Hilbert modular form for $\Gamma_K$, whose divisor $\dv(F)=\sum_{n<0} \tilde c(n)T(-n)$
is given by  Hirzebruch-Zagier divisors. Then there is a weakly holomorphic modular form $f\in W_0^+(D,\chi_D)$
with principal part $\sum_{n<0}c(n)q^n$, and, up to a constant multiple, $F$ is equal to the Borcherds lift of $f$
in the sense of Theorem \ref{hilbert}.
\end{theorem}

As a corollary it can be deduced that the dimension of $\ch^1_{HZ}(\widetilde{X}_K)_\Q$ is equal to the dimension of $M_2^+(D,\chi_D)$ complementing Theorem \ref{hirzebruchzagier}.
It is not hard to see that our map $S_2^+(D,\chi_D)\to S_2(\GK)$ coincides with the celebrated  Doi-Naganuma lift
\cite{DN}, \cite{Na}, \cite{Za}. The above construction can be viewed as a new approach to it. 

The weak converse theorem is proved in much greater generality in \cite{Br2}. Combining the argument of \cite{Br2}
with techniques of \cite{BF} it could probably be extended to hold in full generality for $\Orth(2,\ell)$.
However, for the strong converse theorem the situation seems far more complicated. It is proved in  \cite{Br2} for
modular forms on $\Gamma(L)\subset\Orth(2,\ell)$, where $\Gamma(L)$ denotes the discriminant kernel of the
orthogonal group of an even lattice $L$ of signature $(2,\ell)$ that splits two hyperbolic planes over $\Z$.
For example, if we go to congruence subgroups of the Hilbert modular group $\GK$, it is not clear whether the
analogue of Theorem \ref{converse2} holds or not. 

\section{Automorphic Green functions}\label{sect:greens}

In this section we look at the regularized theta lifts of weak Maass forms from a different perspective. By the
discussion of the previous section, for every positive integer $m$ there exists a unique weak Maass form $f_m\in
\calN^+_0(D,\chi_D)$, whose principal part is equal to $q^{-m}$ if $m\not\equiv 0 \pmod{D}$, and equal to
$\frac{1}{2} q^{-m}$, if $m\equiv 0 \pmod{D}$. The lift
\[
\phi_m(z_1,z_2)=\frac{1}{2}\Phi(z_1,z_2,f_m)
\]
of $f_m$ is a real analytic function on $Y_K$ with a logarithmic singularity along $-2T(m)$.

Here we present a different, more naive, construction of $\phi_m(z_1,z_2)$. For details see \cite{Br1}.
The idea is to construct $\phi_m(z_1,z_2)$ directly as a Poincar\'e series by summing over the logarithms of the
defining equations of $T(m)$. We consider the sum
\begin{align}\label{formalsum}
\sum_{\substack{(a,b,\lambda)\in\Z^2\oplus\frakd^{-1}\\
ab-\lambda\lambda'=m/D}} \log\left|\frac{az_1 \bar{z}_2 +\lambda z_1+\lambda'\bar{z}_2+b}{az_1 z_2 +\lambda
z_1+\lambda'z_2+b}\right|.
\end{align}
The denominators of the summands ensure that this function has a logarithmic singularity along $-2T(m)$ in the same way  as $\phi_m(z_1,z_2)$. The enumerators are smooth on the whole $\H^2$. They are included to make the sum
formally $\GK$-invariant. Unfortunately, the sum diverges. However, it can be regularized in the following way. If
we put $Q_0(z)=\tfrac{1}{2}\log\left(\tfrac{z+1}{z-1}\right)$, we may rewrite  the summands as
\[
\log\left|\frac{az_1 \bar{z}_2 +\lambda z_1+\lambda'\bar{z}_2+b}{az_1 z_2 +\lambda z_1+\lambda'z_2+b}\right| =
Q_{0}\left(1+\frac{ |az_1 z_2 +\lambda z_1+\lambda' z_2+b|^2}{ 2y_1 y_2 m/D }\right).
\]
Now we replace $Q_0$ by the $1$-parameter family $Q_{s-1}$ of Legendre functions of the second kind (cf.~\cite{AS}
\S8), defined by
\begin{equation}\label{legendre}
Q_{s-1}(z)=\int\limits_0^\infty (z+\sqrt{z^2-1} \cosh u)^{-s}du.
\end{equation}
Here $z>1$ and $s\in \C$ with $\Re(s)>0$. If we insert $s=1$, we get back the above $Q_0$. Hence we consider
\begin{equation}\label{phis}
\Phi_m(z_1,z_2,s)  = \sum_{\substack{ a,b\in\Z \\ \lambda\in\frakd^{-1}\\ab- \norm(\lambda)=m/D }}
Q_{s-1}\left(1+\frac{ |az_1 z_2 +\lambda z_1+\lambda' z_2+b|^2}{ 2y_1 y_2 m/D }\right).
\end{equation}
It is easily seen that this series converges normally for $(z_1,z_2)\in\H^2-T(m)$ and $\Re(s)>1$ and therefore
defines a $\GK$-invariant function, which has logarithmic growth along $-2T(m)$. It is an eigenfunction of
the hyperbolic Laplacians $\Delta^{(j)}$ with eigenvalue $s(s-1)$, because of the differential equation satisfied
by $Q_{s-1}$. Notice that for $D=m=1$ the function $\Phi_m(z_1,z_2,s)$ is simply the classical resolvent kernel
for $\Sl_2(\Z)$ (cf.~\cite{Hej}, \cite{Ni}). One can compute the Fourier expansion of $\Phi_m(z_1,z_2,s)$ explicitly and use
it to obtain a meromorphic continuation  to $s\in \C$. At $s=1$ there is a simple pole, reflecting the divergence
of the formal sum \eqref{formalsum}. We define the regularization $\Phi_m(z_1,z_2)$ of \eqref{formalsum} to be the
constant term in the Laurent expansion of $\Phi_m(z_1,z_2,s)$ at $s=1$.

It turns out that $\Phi_m(z_1,z_2)$ is, up to an additive constant $L_m$, equal to the
function $\phi_m(z_1,z_2)$ above (see \cite{Br2} Proposition 2.11 and Theorem 2.14). Here the constant $L_m$
\label{lm} is quite interesting, since it is given by the derivative of the $m$-th coefficient of a certain
Eisenstein series $E_2(\tau,s)$ of weight $2$ for $\Gamma_0(D)$ and $\chi_D$ \cite{BrKue}, \cite{BBK}. We will come back to this in section
\ref{intgreen}.

One may use the Fourier expansion of
$\Phi_m(z_1,z_2)$ and identities for certain finite exponential sums of \cite{Za} to obtain a different independent
proof of Theorem \ref{hilbert}.

The following integral formula is fundamental (see \cite{BrKue} Theorem 4.7, \cite{BBK}). It justifies why
$\Phi_m(z_1,z_2,s)$ (and also $\Phi_m(z_1,z_2)$) is called an {\em automorphic Green function} for the divisor $T(m)$.

\begin{theorem}\label{intphif}
Let $h:Y_K\to\C$ be a bounded eigenfunction of the Laplacian $\Delta^{(1)}$ (or $\Delta^{(2)}$) with
  eigenvalue $\lambda$.  Then for $s\in \C$ with $\Re(s)>1$ we have
\[
\int\limits_{\GK\bs \H^2} \Phi_m(z_1,z_2,s) h(z_1,z_2)\, \Omega^2 = \frac{1}{s(s-1)-\lambda} \int\limits_{T(m)}
h(z_1,z_2) \,\Omega.
\]
Here $\Omega$ is defined by \eqref{defomega} so that $\Omega^2$ is an invariant volume form on $Y_K$.
\end{theorem}


Such automorphic Green functions are constructed in greater generality for $\Orth(2,\ell)$ in \cite{Br2} using
the regularized theta lift of Hejhal-Poincar\'e series, and independently in \cite{OT} from the point of view of
spherical functions on real Lie groups.

\section{Integrals of automorphic Green functions}\label{intgreen}
\label{sect:61}
It is well known that the volume of $Y_K$ is given by $\vol_\Omega(Y_K)= \int_{Y_K} \Omega^2=\zeta_K(-1)$, where
$\zeta_K(s)$ denotes the Dedekind zeta function of $K$. The volume of a divisor $C$ on $Y_K$ is defined as the
integral $\vol_\Omega(C)= \int_{C} \Omega$. One can show that the integral is finite, see e.g.~\cite{Br3}. If $C$ is effective, then its volume is
positive. 

It is a well known fact that the volumes of Hirzebruch-Zagier divisors are given by the Fourier coefficients of
the unique normalized Eisenstein series in $M_2^+(D,\chi_D)$ (see \cite{Fra}, \cite{Ha}, and \cite{Ge} chapter V.5). Let
us recall the definition of that Eisenstein series. In weight $k$ there are the two non-holomorphic Eisenstein
series
\begin{align*}
E_k^\infty(\tau,s)&= \sum_{\substack{c,d\in \Z\\c\equiv 0 \;(D)}} \chi_D(d) \frac{1}{(c\tau + d)^k} \frac{y^s} {|c\tau + d|^{2s}},\\
E_k^0(\tau,s)&= \sum_{c,d\in \Z} \chi_D(c) \frac{1}{(c\tau + d)^k} \frac{y^s} {|c\tau + d|^{2s}}
\end{align*}
for $\Gamma_0(D)$ with character $\chi_D$, the former corresponding to the cusp $\infty$ of $\Gamma_0(D)$, the
latter to the cusp $0$. (By our assumption that $D$ be prime, these are the only cusps of $\Gamma_0(D)$.) They
converge for $\Re(s)>1-k/2$ and have a meromorphic continuation in $s$ to the full complex plane. If $k\geq 2$,
the special values $E_k^\infty(\tau,0)$ and $E_k^0(\tau,0)$ are holomorphic in $\tau$ and define elements of
$M_k(D,\chi_D)$. One can show that the linear combination
\[
E_k(\tau,s)=\frac{1}{2 L(k+2s,\chi_D)}\left( D^s E_k^\infty(\tau,s)+ D^{1/2-k-s}E_k^0(\tau,s) \right)
\]
satisfies the plus space condition. (This follows most easily from Lemma 3 of \cite{BB}.) Here $L(s,\chi_D)$
denotes the $L$-series associated with the Dirichlet character $\chi_D$. In particular we have $E_k(\tau,0)\in
M_k^+(D,\chi_D)$. The Fourier expansion of $E_k(\tau,s)$ has the form
\begin{align}
\label{eiss} E_k(\tau,s) = \sum_{n\in\Z} C(n,s)\calW_s(4\pi n v) e^{2\pi i nu},
\end{align}
where the $C(n,s)$ are complex coefficients independent of $v$; and $\calW_s(v)$ is a certain Whittaker function,
which we normalize as in \cite{BrKue} (3.2). The precise normalization is not important for our purposes here, we
only need that it is a universal function for all $n$ of the same sign. The coefficients $C(n,s)$ are computed
for instance in \cite{BrKue} section~5, Example~2. Here we only state the special value
\begin{equation}\label{eis}
E_k(\tau,0) = 1+\sum_{n\geq 1} C(n,0)q^n=  1+\frac{2}{L(1-k,\chi_D)} \sum_{n\geq 1} \sum_{d\mid n} d^{k-1} \left(
\chi_D(d) + \chi_D(n/d)\right) q^n,
\end{equation}
which is obtained in the standard way (see \cite{He}, Werke p.~818) using the functional equation of $L(s,\chi_D)$.

\begin{theorem}\label{hzvol}
We have
\[
E_2(\tau,0)= 1 - \frac{2}{\vol_\Omega(Y_K)} \sum_{m\geq 1} \vol_\Omega(T(m)) q^m.
\]
\end{theorem}
Similar identities hold in much greater generality for special cycles on arithmetic quotients of $\Orth(p,q)$ and
$\Uni(p,q)$,  see e.g.~\cite{Ge2}, \cite{Ku:Borcherds},  \cite{Ku:Bourbaki}, \cite{Oda:SIEGEL}. (Observe that our normalization of $\vol_\Omega(T(m))$ equals twice the volume of $T(m)$ in \cite{BBK}.)

Let us briefly indicate, how Theorem \ref{hzvol} can be deduced from the properties of the automorphic Green
functions $\Phi_m(z_1,z_2,s)$. For instance, from the description
as a regularized theta lift it follows that the residue at $s=1$ of $\Phi_m(z_1,z_2,s)$ is equal to the constant
coefficient $a_m(0)$ of the weak Maass form $f_m\in \calN^+_0(D,\chi_D)$ defined at the beginning of section
\ref{sect:greens}. By means of the relationship of the spaces $\calN^+_k(D,\chi_D)$ and $M_{2-k}^+(D,\chi_D)$,
which is also implicit in \eqref{exseq}, one finds that $a_m(0)=-\frac{1}{2}C(m,0)$ (see \cite{BF} Proposition
3.5). Therefore we have
\begin{equation*}
\Phi_m(z_1,z_2)=\lim_{s\to 1}\left( \Phi_m(z_1,z_2,s)+\frac{C(m,0)}{2(s-1)}\right).
\end{equation*}
Using growth estimates for $\Phi_m(z_1,z_2,s)$, which can be deduced from the constant coefficients of the Fourier
expansions, we obtain:

\begin{proposition}\label{intphi}
The function $\Phi_m(z_1,z_2)$ belongs to $L^p(Y_K,\Omega^2)$ for any $p<2$, and 
\[
\int_{Y_K} \Phi_m(z_1,z_2) \,\Omega^2 =\lim_{s\to 1}\int_{Y_K} \left( \Phi_m(z_1,z_2,s)+\frac{C(m,0)}{2(s-1)}\right) \Omega^2.
\]
\end{proposition}

If we apply Theorem \ref{intphif} for the constant function $h=1$, we may compute the integral. It is equal to
\begin{align}\label{h1}
\frac{1}{s(s-1)} \int_{T(m)} \Omega + \frac{C(m,0)}{2(s-1)} \int_{Y_K} \Omega^2 = \frac{1}{(s-1)} \left(
\frac{\vol_\Omega(T(m))}{s} + \frac{C(m,0)}{2}\vol_\Omega(Y_K)\right).
\end{align}
Since the limit $s\to 1$ exists, the quantity in parenthesis on the right hand side has to vanish at $s=1$. This
yields the assertion of Theorem \ref{hzvol}.

So far we have essentially exploited the existence of the integral $\int_{Y_K} \Phi_m(z_1,z_2) \,\Omega^2$, which
means that the residue in the Laurent expansion of \eqref{h1} at $s=1$ vanishes. We may actually compute the
constant term of that expansion, that is, the value of the integral. It is equal to $-\vol_\Omega(T(m))$.

One can further improve this result by observing that the full coefficient $C(m,s)$ as a function of $s$  occurs
in the constant term of $\Phi_m(z_1,z_2,s)$. More precisely, if we define
\begin{align}
\label{defgm} G_m(z_1,z_2)= \frac{1}{2}\lim_{s\to 1} \big( \Phi_m(z_1,z_2,s)+B(s)\zeta(2s-1) C(m,s-1)\big),
\end{align}
with
\[
B(s)=  \frac{(16\pi)^{s-1}\Gamma(s-1/2) s}{\Gamma(1/2)(2s-1)},
\]
then one can show that $G_m(z_1,z_2)=\frac{1}{4}\Phi(z_1,z_2,f_m)+\frac{1}{2}a_m(0)(\Gamma'(1)-\log(8\pi))$ (which is
essentially the calculation of the constant $L_m$ on page \pageref{lm}). This means in particular that if $F$ is
the Borcherds lift of a weakly holomorphic modular form $f$ with coefficients $a(n)$, then its {\em Petersson
metric}
is given by
\begin{align} \label{krolim}
\log\|F(z_1,z_2)\|_{Pet} := \log \left(|F(z_1,z_2)| (16\pi^2 y_1 y_2)^{a(0)/2}\right)=- \sum_{n<0} \tilde a(n) G_{-n}(z_1,z_2).
\end{align}
The latter identity  can be viewed as a generalization of the Kronecker limit formula expressing the logarithm
of the absolute value of the discriminant function \eqref{defdelta} as the constant term in the Laurent expansion
at $s=1$ of the non-holomorphic Eisenstein series of weight $0$ for $\Sl_2(\Z)$ (see \cite{BrKue} (4.14)). Notice
that the constant $B(s)$ in \eqref{defgm} does not depend on $m$. It changes if the normalization of the Whittaker
function $\calW_s(v)$ is varied.  Arguing as above we find that (see \cite{BrKue} Theorem 4.10)
\begin{align*}
\int_{Y_K} G_m(z_1,z_2) \,\Omega^2 
&= -\frac{\vol_\Omega(T(m))}{2}\big(C'(m,0)/C(m,0)+\log(4\pi) -\Gamma'(1)\big).
\end{align*}
If we insert the explicit formula for $C(m,s)$, we get
\begin{align}\label{greenint}
\int_{Y_K} G_m(z_1,z_2) \,\Omega^2
&=-\vol_\Omega(T(m)) \left( \frac{L'(-1,\chi_D)}{L(-1,\chi_D)}+\frac{1}{2}-\frac{\sigma_m'(-1)}{\sigma_m(-1)}+\frac{1}{2}\log(D)\right),
\end{align}
where 
\begin{align}\label{def:sig}
\sigma_m(s)=m^{(1-s)/2}\sum_{d\mid m} d^{s} \left( \chi_D(d) + \chi_D(m/d)\right). 
\end{align}

In particular, in view of
\eqref{krolim}, the integral over the logarithm of the Petersson metric of any Borcherds product can be computed
explicitly (see also \cite{Ku:Borcherds}). For example, if $K=\Q(\sqrt{5})$, we obtain for the Gundlach theta function
\[
\int_{Y_K} \log\left( 2^{-6}|\Theta(z_1,z_2)| (16\pi^2 y_1 y_2)^{5/2}\right) \Omega^2
=-\zeta(-1)\left(2\frac{L'(-1,\chi_D)}{L(-1,\chi_D)}+1+\log(5)\right).
\]
Such integrals play a fundamental role in the Arakelov intersection theory of Hirzebruch-Zagier divisors. We will
come back to that in the section \ref{sect:8}.

The integral of the logarithm of the Petersson metric of a Borcherds product was first calculated by Kudla in
\cite{Ku:Borcherds} using a different approach based on the Siegel-Weil formula. We recall that the quantity
$\Phi(z_1,z_2,f)$ we want to integrate is given by the theta integral \eqref{regint} of a weakly holomorphic
modular form $f$. Now the idea is to interchange the $(z_1,z_2)$-integration with the regularized integration over
$\tau$ and to compute
\begin{align}\label{kudlaint}
\int_\calF \left\langle f(\tau),\int_{Y_K} \Theta_S(\tau;z_1,z_2) \,d\mu_z \right\rangle\,d\mu_\tau,
\end{align}
where $d\mu_z$ denotes the invariant measure on $Y_K$.  (Notice that this needs a careful justification.) The inner
integral over the Siegel theta function can be determined by means of the Siegel-Weil formula. It yields an
Eisenstein series of weight $0$ for $\Gamma_0(D)$, which can be written in terms of the lowering operator and our
Eisenstein series \eqref{eiss} of weight $2$ as $\frac{1}{s-1}L_2 E_2(\tau,s)$. The integrand for the remaining
regularized integral over $\tau$ is now essentially $\frac{1}{s-1} d(\langle f (\tau), E_2(\tau,s)\rangle d\tau) $ at
$s=1$, so that we may use Stoke's theorem to compute it. The derivative of $E_2(\tau,s)$ occurs because of the
factor $\frac{1}{s-1}$.

We conclude this section by giving a characterization of the automorphic Green function $\Phi_m(z_1,z_2)$. It can be proved in a similar way as Theorem \ref{converse1}.

\begin{proposition}
Let $f$ be a smooth function on $Y_K- T(m)$ with the properties:
\begin{enumerate}
\item[(i)] $f$ has a logarithmic singularity along $T(m)$, 
\item[(ii)] $(\Delta^{(1)}+\Delta^{(2)}) f = \text{constant}$,
\item[(iii)] $f\in L^{1+\eps}(Y_K,\Omega^2)$ for some $\eps>0$,
\item[(iv)] $\int_{Y_K} f(z_1,z_2)\, \Omega^2 = \frac{1}{2}\vol_\Omega (T(m))$.
\end{enumerate}
Then $f(z_1,z_2)=-\frac{1}{2}\Phi_m(z_1,z_2)$.
\end{proposition}

\section{Arithmetic of Hirzebruch-Zagier divisors}
\label{sect:7}


In their paper on the intersection of modular correspondences, Gross and Keating interpreted classical results of Hurwitz and Kronecker by the observation that the intersection number of two modular correspondences on $Y_{\Q\oplus\Q}=\Gamma(1)^2 \bs \H^2$ is given by the coefficients of the classical Siegel Eisenstein series $E^{(2)}(Z,s)$ of weight 2 and genus 2 at $s=0$. Their main result was that the arithmetic intersection numbers of three such modular correspondences on the regular model $\Spec \Z[j,j']$ of $Y_{\Q\oplus\Q}$ is given by the coefficients of the derivative of the Siegel Eisenstein series $E^{(3)}(Z,s)$ of weight 2 and genus 3 at $s=0$ \cite{GK}.
Observe that $Y_{\Q\oplus\Q}$ can be viewed as the ``degenerate'' Hilbert modular surface with discriminant $D=1$ and the modular correspondences as Hirzebruch-Zagier divisors in this case. 

Kudla proved that the arithmetic intersection numbers 
in the sense of Arakelov geometry 
of certain arithmetic special divisors on a regular model of a Shimura curve are dictated by the coefficients of the derivative of a Siegel Eisenstein series of weight 2 and genus 2 at $s=0$ \cite{Ku1}. (So Kudla considers an arithmetic surfaces, rather than an arithmetic $3$-fold as in the case of Gross and Keating. This explains the different genus.)
Here the arithmetic divisors are pairs consisting of a special divisor on the regular model and a certain Green function for the induced divisor on the corresponding complex variety, fitting in the setup of arithmetic intersection theory as in \cite{SABK}.

In further works Kudla, Rapoport, and Yang developed an extensive program relating arithmetic special divisors on Shimura varieties of type $\Orth(2,\ell)$ and their arithmetic intersection theory to automorphic forms, in particular to the coefficients of the derivatives of Siegel Eisenstein series. Most of this is conjectural, but in important special cases these conjectures are meanwhile proved. (See e.g.~\cite{Ku:Harvard} for the $\Orth(2,1)$ case of Shimura curves, \cite{KRY} for the $\Orth(2,0)$ case of CM elliptic curves, \cite{KR} for partial results in the $\Orth(2,3)$ case of Siegel modular threefolds, and \cite{Ku:MSRI} for an overview.)
Notice that $Y_{\Q\oplus\Q}$ can be described in terms of $\Orth(2,2)$.

One conclusion of this general picture is that the geometric results over $\C$ of Hirzebruch and Zagier (as e.g.~Theorems \ref{hirzebruchzagier} and \ref{hzvol}) and their generalizations to $\Orth(2,\ell)$ should have arithmetic analogues over $\Z$. Here the classical intersection theory has to be replaced by Arakelov intersection theory.

In this section we discuss, how Borcherds products can be used to obtain new results in that direction. We begin by recalling some facts on the arithmetic of Hilbert modular surfaces. 

\medskip

In section \ref{sect:2} we briefly discussed that modular curves have a moduli interpretation as a moduli space for isomorphism classes of elliptic curves with additional structure. The same is true for Hilbert modular surfaces, which is the starting point for arithmetic investigations.
More precisely, $Y_K$ parametrizes isomorphism classes of triples $(A, \iota,\psi)$, where $A$ is an abelian surface over $\C$, $\iota$ is an $\OK$-multiplication, that is, a ring homomorphism $\OK\to \End(A)$, and $\psi$ is a $\frakd^{-1}$-polarization, that is, an isomorphism of $\OK$-modules $\frakd^{-1}\to
 \Hom_{\OK}(A,A^\vee)^{\sym} $ 
from the inverse different $\frakd^{-1}=\frac{1}{\sqrt{D}}\OK$ to the module of $\OK$-linear symmetric homomorphisms, taking the totally positive elements of  $\frakd^{-1}$ to  $\OK$-linear polarizations (see \cite{Go} Chapter 2).  

The moduli description now makes sense over any scheme $S$ over $\Z$. (Here one has to require that $\psi$ fulfill an extra technical condition called the Deligne-Pappas condition, see \cite{DePa}. That condition is automatically fulfilled in characteristic $0$.) 
Due to the work of Rapoport, Deligne, and Pappas it is known that
the moduli problem ``Abelian surfaces over $S$ with $\OK$-multiplication
and $\frakd^{-1}$-polarization with Deligne-Pappas condition'' 
is represented by a regular
algebraic stack $\mathcal{H}$, which is flat and of
relative dimension two over $\Spec \Z$. It is smooth
over $\Spec \Z[1/D]$, and the fiber of
$\mathcal{H}$ above $D$ is smooth outside a closed subset
of codimension $2$.

The corresponding complex variety $\calH(\C)$ is isomorphic to $Y_K$. The isomorphism is obtained by associating to $z=(z_1,z_2) \in \H^2$
the  abelian surface $A_z=\C^2/\Lambda_{z}$ over $\C$ given by the lattice
\[
\Lambda_z 
= \left\{
\begin{pmatrix} \alpha z_1 +  \beta \\
\alpha' z_2 +  \beta'
\end{pmatrix} 
\in \C^2;\quad \alpha,\beta\in \OK\right\}\subset \C^2,
\]
together with the $\OK$-multiplication $\iota$ induced by the natural action $\iota(\nu)=\kzxz{\nu}{0}{0}{\nu'}$ of $\OK$ on $\C^2$, and a certain $\frakd^{-1}$-polarization.

For $k\in \Z$ sufficiently divisible there exists a line bundle $\calM_k$ on $\calH$ (the $k$-th
power of the pull-back along the zero section 
of the determinant of the relative cotangent
bundle of the universal family over $\mathcal{H}$) such that the induced bundle on $\calH(\C)$ can be identified with the line bundle $\calM_k(\C)$ of Hilbert modular forms of weight $k$ for $\GK$ of the previous sections. By the $q$-expansion principle and the Koecher principle, the global sections of $\calM_k$ can be identified with Hilbert modular forms of weight $k$ for $\GK$ with integral rational Fourier coefficients.
 
There exists an arithmetic Baily-Borel compactification $\overline{\calH}$ of the coarse moduli space corresponding to $\calH$, which can be described as
\begin{align}\label{eq:minhilbert}
\overline{\calH} = 
\Proj   \bigg(\bigoplus_{k} 
H^0 ( \calH,  \calM_k)\bigg).
\end{align}
The scheme $\overline{\calH}$ is  normal, projective,
and flat over $\Spec\Z$ (see \cite{Ch}, p.~549), and $\overline{\calH}(\C)\cong X_K$. 
Furthermore, its fibers over $\Spec \Z$ are 
irreducible (see \cite{DePa}, p.~65). 
By construction, the bundle $\calM_k$ extends to $\overline \calH$. 

Throughout the rest of this paper we will make the following  

\begin{assumption}\label{ass:model}
There exists a desingularization $\pi:\widetilde{\calX}_K \to \overline{\calH}$ 
by a regular scheme $\widetilde{\calX}_K$, which is projective and flat over $\Z$, 
such that 
the regular locus $\overline{\calH}^{reg}$ is fiber-wise dense in $\widetilde{\calX}_K$, and such that the induced morphism $\widetilde{\calX}_K (\C)\to X_K$ is a desingularization as in the previous sections $\widetilde{X}_K$.
\end{assumption}

This assumption simplifies the exposition (it might actually be too optimistic). 
Notice that the singularities of $\overline{\calH}$ at the boundary can be resolved by considering a suitable toroidal compactification of $\calH$. So only the singularities corresponding to elliptic fixed points would need to be resolved. 
If one wants to obtain unconditional results one can impose an additional level structure in order to get a fine moduli problem and work with a suitable toroidal compactification of the corresponding moduli scheme (as is done in \cite{BBK}). Unfortunately, in that way one only gets a regular scheme $\widetilde{\calH}(N)$ which is projective and flat over $\Z[\zeta_N, 1/N]$, where $\zeta_N$ denotes a primitive $N$-th root of unity and $N\geq 3$ the level. 

We define the line bundle of modular forms of weight $k$ on 
$\widetilde{\calX}_K$ as the pullback $\pi^*(\calM_k)$. For simplicity we will also denote it by $\calM_k$.

It can be shown that the Hirzebruch-Zagier divisors on $X_K$ are defined over $\Q$, that is, $T(m)$ is obtained by base change from a divisor on the generic fiber $\overline\calH\times_\Z\Q$ of $\overline{\calH}$. 
We define the Hirzebruch-Zagier divisor $T(m)$ on the generic fiber $\widetilde{\calX}_K\times_\Z\Q$ as the pullback of $T(m)$ on $\overline\calH\times_\Z\Q$.
Moreover, we define the
Hirzebruch-Zagier divisor $\mathcal{T}(m)$ on
$\widetilde{\calX}_K$ as the Zariski closure of $T(m)$.

\medskip

We now briefly recall some basic properties of arithmetic Chow rings (see e.g.~\cite{SABK}). Since $Y_K$ is non-compact, the natural metrics on  automorphic vector bundles have singularities at the boundary \cite{Mu}, \cite{BKK2}. Therefore we need to work with the extended arithmetic Chow ring $\widehat{\textrm{CH}}^*(\mathcal{X},\mathcal{D}_{\pre})$ constructed in \cite{BKK}. In this ring the Green objects satisfy beside the usual logarithmic additional log-log growth conditions.

Let $\calX$ be an arithmetic variety over $\Z$, i.e., a regular scheme, which is projective and flat over $\Z$. Moreover, let $\D$ be a fixed normal crossing divisor on the complex variety $\calX(\C)$, which is stable under complex conjugation.
An arithmetic divisor on $\calX$ (in the sense of \cite{BKK}) is a pair 
\[
(y, \frakg_y),
\] 
where $y$ is a divisor on the scheme $\calX$ and $\frakg_y$ is a pre-log-log Green object for  the induced divisor $y(\C)$ on $\calX(\C)$.
In particular, a pair $(y,g_y)$ where $g_y$ is a {\em pre-log-log Green function} for $y$, determines an arithmetic divisor. This essentially means that  $g_y$ is a smooth function on $\calX(\C)-(y(\C)\cup \D)$, invariant under complex conjugation, with logarithmic singularities along the irreducible components of $y(\C)$ and pre-log-log singularities along $\D$ such that the $\partial\bar\partial$-equation of currents holds:
\[
-2\partial\bar\partial [g] = [-2\partial\bar\partial g]-\delta_y.
\]
Here $[\cdot]$ denotes the current associated to a differential form and $\delta_y$ the Dirac current for $y$  normalized as in \cite{BKK} and \cite{BBK}. A differential  form $\alpha$ is called pre-log singular (pre-log-log singular), if $\alpha$, $\partial\alpha$, $\bar\partial\alpha$, and $\partial\bar\partial \alpha$ have only logarithmic growth (respectively log-log growth). 
We write $\widehat{\rm Z}^1(\calX,\calD_\pre)$ for the free abelian group generated by the arithmetic divisors on $\calX$.
Here $\calD_\pre$ stands for the Deligne algebra with pre-log-log forms along $\D$, which is needed for the precise description of Green objects in \cite{BKK}.
Moreover, we write $\rata^1(\calX)$ for the subgroup of $\za^1(\calX,\mathcal{D}_{\pre})$ given by arithmetic divisors of the form $(\dv(f), -\log|f|)$, where $f$ is a rational function on $\calX$ and $|f|$ the absolute value of the induced function on $\calX(\C)$.
The first {\em arithmetic Chow group} of $\calX$ with
log-log growth along $\D$ is defined by
\[
\cha^1(\calX,\mathcal{D}_{\pre})=\za^1(\calX,\mathcal{D}_{\pre}) 
 \big/ \rata^1(\calX). 
\]

More generally, in \cite{BKK} arithmetic Chow groups $\cha^p(\calX,\mathcal{D}_{\pre})$ of codimension $p$ arithmetic cycles with log-log growth along $\D$ are defined. 
There exists an arithmetic intersection  product 
$$
\cha^p(\calX, \mathcal{D}_{\pre}) \otimes \cha^q(\calX,
\mathcal{D}_{\pre} ) \longrightarrow \cha^{p+q}(\calX,
\mathcal{D}_{\pre} )_\Q,
$$
and 
$$\cha^*(\calX, \mathcal{D}_{\pre})_{\mathbb{Q}}=\bigoplus_{p\ge 0}
\cha^p(\calX, \mathcal{D}_{\pre})\otimes_\Z\mathbb{Q}
$$ 
equipped with this product
has the structure of a commutative associative ring. 

For instance, if $\calX=\Spec \Z$, then the closed points of $\Spec \Z$ can be identified with the primes of $\Z$. An arithmetic divisor is a pair $(\sum_p n_p p,\, g)$ consisting of a finite formal $\Z$-linear combination of primes $p$ and a real number $g$. The elements of  $\rata^1(\Spec \Z)$ are the pairs of the form $(\sum_{p\mid N} \ord_p(N) p,\, -\log|N| )$ for $N\in \Q$.
This implies that 
\[
\dega: \cha^1(\Spec \Z, \mathcal{D}_{\pre})\longrightarrow \R, \quad \left(\sum n_p p, \,g \right) \mapsto g +\sum n_p \log(p)
\]
is an isomorphism. It is common to identify $\cha^1(\Spec \Z, \mathcal{D}_{\pre})$ with $\R$.

There also is an arithmetic analogue of the Picard group: The {\em arithmetic Picard group} $\pica(\calX, \mathcal{D}_{\pre})$ is the group of isomorphism classes of pre-log singular hermitian line bundles on $\calX$. Here a pre-log singular hermitian line bundle is a pair $\overline{\calL}=(\calL, \|\cdot\|)$ consisting of a line bundle $\calL$ on $\calX$, and a smooth hermitian metric $\|\cdot\|$ on the induced complex line bundle on $\calX(\C)-\D$, invariant under complex conjugation, and such that $-\log\|s\|$ has logarithmic singularities along $\dv(s)(\C)$ and pre-log-log singularities along $\D$ for any rational section $s$ of $\calL$. 
If $\overline{\calL}$ is a pre-log singular hermitian line bundle and $s$ a rational section of $\calL$, then, essentially by the Poincar\'e-Lelong 
lemma, 
\[
\cca_1(\overline{\calL})=(\dv(s), -\log\|s\|) 
\]
defines a class in $\cha^1(\calX, \mathcal{D}_{\pre})$, which is independent of the choice of $s$. It is called the first arithmetic Chern class of $\overline{\calL}$. The assignment $\overline{\calL}\mapsto\cca_1(\overline{\calL})$ actually induces an isomorphism 
\begin{align*}
 \cca_1: \pica(\calX,\mathcal{D}_\pre) \longrightarrow
  \cha^1(\calX,\mathcal{D}_\pre).
\end{align*}

\medskip

We now consider the arithmetic Chow ring $\cha^*(\widetilde{\calX}_K, \mathcal{D}_{\pre})_{\mathbb{Q}}$ of the model $\widetilde{\calX}_K$ of our Hilbert modular surface, where we take for $\D$ the normal crossing divisor $\D_K=\pi^{-1}(X_K^{sing})$. For details we refer to \cite{BBK}.
The Green functions of section \ref{sect:greens} turn out to be particularly nice, because they fit into the arithmetic Chow theory of $\widetilde{\calX}_K$.

\begin{theorem} \label{thm:defcycle} 
The pair 
\[
\widehat{\calT}(m)=\left(\calT(m), G_m\right)
\]
defines an element of $\cha^1(\widetilde{\calX}_K, \mathcal{D}_{\pre})$, called the arithmetic Hirzebruch-Zagier divisor of discriminant $m$. Here $G_m$ is the automorphic Green function defined by \eqref{defgm}.
\end{theorem}

Notice that $G_m$ always has log-log singularities along $\D_K$, even if $T(m)$ is disjoint to $\D_K$. So $G_m$ does not define a Green function for $T(m)$ in the classical arithmetic Chow theory due to Gillet and Soul\'e. We therefore really need the extension of \cite{BKK}. Observe that the arithmetic divisors of Theorem \ref{thm:defcycle} slightly differ from those considered by Kudla, Rapoport, and Yang. For instance, they often contain boundary components (which is possible since $Y_K$ is non-compact), and are built with different Green functions.

Moreover, we obtain an element of  $\pica(\widetilde{\calX}_K, \mathcal{D}_{\pre})$, by equipping 
the line bundle of modular forms with the {\em Petersson metric}. Recall that if $F\in \calM_k(\C)(U)$ is a rational section over an open subset $U\subset Y_K$, then its Petersson metric is given by
\begin{align*}
\|F(z_1,z_2)\|^2_{\Pet} 
=|F(z_1,z_2)|^2 (16 \pi^2 y_1 y_2)^k.
\end{align*}
This defines a pre-log singular hermitian metric on $\calM_k(\C)$ (with  respect to  $\D_K$). 
We denote  the corresponding pre-log singular hermitian line bundle by $\overline{\calM}_k=(\calM_k,\|\cdot\|^2_{\Pet})$.
(That the Petersson metric has singularities at the boundary is easily seen: For instance, if $z=(z_1,z_2)$ approaches the cusp $\infty$ of $Y_K$, then $y_1 y_2 \to \infty$ by construction of the Baily-Borel topology. At the elliptic fixed points it is continuous, but the derivatives do have singularities.)

A central idea in \cite{BBK} is to connect the arithmetic of Borcherds products and the properties of the automorphic Green functions $G_m$ to derive information on $\cha_{HZ}^1(\widetilde{\calX}_K, \mathcal{D}_{\pre})_\Q$, the subspace of $\cha^1(\widetilde{\calX}_K, \mathcal{D}_{\pre})_\Q$ spanned by the arithmetic Hirzebruch-Zagier divisors.

\begin{theorem} 
\label{thm:mor}
Recall Assumption \ref{ass:model}.
The homomorphism $\z^1_{HZ}(\widetilde{X}_K) \to \za^1_{HZ}(\widetilde{\calX}_K,\mathcal{D}_{\pre})$ defined by $T(m)\mapsto \widehat{\calT}(m)$ induces an isomorphism
\[
 \ch_{HZ}^1(\widetilde{X}_K)_\Q \longrightarrow
 \cha_{HZ}^1(\widetilde{\calX}_K, \mathcal{D}_{\pre})_\Q,
\]
taking $\cc_1(\calM_k(\C))$ to $\cca_1(\calM_k)$.   
\end{theorem}

\begin{proof}[Sketch of the proof.]
We have to show that if there is a relation in $\z^1(\widetilde{X}_K)$ among the $T(m)$, we can lift it to a relation in $\za^1(\widetilde{\calX}_K,\mathcal{D}_{\pre})$, and that every relation among arithmetic Hirzebruch-Zagier divisors arises in that way (up to torsion).

So suppose that $F$ is a rational function on $\widetilde{X}_K$ with divisor $\sum_{n<0} \tilde{c}(n) T(-n)$. Then by the strong converse theorem (Theorem \ref{converse2}), we may assume that $F$ is a Borcherds product, that is, the lift of a weakly holomorphic modular form $f\in W_0^+(D,\chi_D)$ with Fourier expansion $\sum_{n} c(n)q^n$ as in Theorem \ref{hilbert}. 
It can be shown that any meromorphic Borcherds product is the quotient of two holomorphic ones (\cite{BBK} Proposition 4.5). Therefore we may write $F=F_1/F_2$, where $F_1$, $F_2$ are holomorphic Borcherds products of the same weight. 
But then the Borcherds product expansion (Theorem \ref{hilbert} (iii)) implies that a positive power of $F_j$ has integral rational Fourier coefficients. Without loss of generality we may assume that already the $F_j$ 
have integral rational Fourier coefficients.
According to the $q$-expansion principle $F_j$ defines a section $\calF_j$ of $\calM_{k}$ on the model $\widetilde{\calX}_K$. Hence the quotient of these sections is a rational function on $\widetilde{\calX}_K$ that specializes to $F$ on the generic fiber.

We claim that the divisor of $\calF_j$ on $\widetilde{\calX}_K$ is horizontal. To see this, we notice that by work of Rapoport, Deligne, and Pappas, the geometric fibers of $\calH$ are irreducible (see \cite{Ra}, \cite{DePa}). It follows by Assumption \ref{ass:model}, that the same holds for the geometric fibers of $\widetilde{\calX}_K$.
Suppose that $\dv(\calF_j)$ contains a vertical component above a prime $p$. 
Then, because of the irreducibility of the fibers, $\dv(\calF_j)$ contains the full fiber above $p$. By the $q$-expansion principle, this implies that all Fourier coefficients of $F_j$ are divisible by $p$. But the Borcherds product expansion of $F_j$ shows that the coefficients are coprime (in fact, the coefficient corresponding to the Weyl vector $\rho(F_j)$ is $1$), and therefore a contradiction.

Thus the divisor of $\calF_1/\calF_2$ is horizontal and equal to $\sum_{n<0} \tilde{c}(n) \calT(-n)$. In view of \eqref{krolim} we may conclude that
\[
\sum_{n<0} \tilde{c}(n) \widehat{\calT}(-n)= \left(\dv(\calF_1/\calF_2), -\log| F|\right) =0\in  \cha^1(\widetilde{\calX}_K,\mathcal{D}_{\pre})_\Q.
\]

Conversely, 
every relation among arithmetic Hirzebruch-Zagier divisors obviously specializes to a relation on the generic fiber.
\end{proof}

As a corollary we see that $\dim(\cha_{HZ}^1(\widetilde{\calX}_K, \mathcal{D}_{\pre})_\Q)=\dim(\ch_{HZ}^1(\widetilde{X}_K)_\Q) =[\frac{D+19}{24}]$. Moreover, diagram \eqref{eq:diag} has an arithmetic analogue, where one has to replace $\z^1_{HZ}(\widetilde{X}_K)$ by  $\za^1_{HZ}(\widetilde{\calX}_K,\mathcal{D}_{\pre})$ and $\ch_{HZ}^1(\widetilde{X}_K)_\Q$ by $\cha_{HZ}^1(\widetilde{\calX}_K, \mathcal{D}_{\pre})_\Q$. Finally, in view of Theorem \ref{hirzebruchzagier}, one 
obtains the following arithmetic Hirzebruch-Zagier theorem
(cf.~\cite{BBK} Theorem 6.2):

\begin{theorem}\label{thm:ags}
The arithmetic generating series
\begin{align} \label{agenseries}
  \widehat{A}(\tau)=\cca_1(\overline{\mathcal{M}}_{-1/2}) +
  \sum_{m>0} \widehat{\mathcal{T}}(m) q^m
\end{align}
is a holomorphic modular form in $M^+_2(D,\chi_D)$ with values in
$\cha^1(\widetilde{\calX}_K,\mathcal{D}_\pre)_\Q$, i.e., an element
of
$M^+_2(D,\chi_D)\otimes_\Q\cha^1(\widetilde{\calX}_K,\mathcal{D}_\pre)_\Q$.
\end{theorem}


\section{Arithmetic intersection numbers}
\label{sect:8}

The first Chern form of the line bundle $\calM_k(\C)$ equipped with the Petersson metric is equal to
\[
\cc_1(\calM_k(\C), \|\cdot\|_\Pet) = 2\pi i k\cdot \Omega,
\]
where $\Omega$ denotes the  K\"ahler form \eqref{defomega}.
Consequently, $\vol_\Omega(Y_K)=\zeta_K(-1)$ can also be regarded as the geometric self intersection number $\calM_1(\C)^2$ of the line bundle of modular forms of weight $1$. Moreover, 
Theorem \ref{hzvol} can be rephrased by saying that the intersection
of the geometric generating series \eqref{genseries} and $\cc_1(\calM_k(\C))$ is given by
\begin{align*}
A(\tau)\cdot \cc_1(\calM_k(\C)) 
&= -\frac{k}{2} \zeta_K(-1)  \cdot E_2(\tau,0),
\end{align*}
where $E_2(\tau,0)\in M_2^+(D,\chi_D)$ is the Eisenstein series \eqref{eis}.

In view of this result it is natural to ask, what the intersection of the arithmetic generating series \eqref{agenseries} with the class $\cca_1(\overline{\calM}_k)^2\in \cha^2(\widetilde{\calX}_K,\mathcal{D}_\pre)_\Q$ is.

\begin{theorem}
\label{thm:B}
Recall Assumption \ref{ass:model}.
We have the following identities of arithmetic intersection numbers:
\begin{align}\label{eq:ags}
\widehat{A}(\tau) \cdot 
\cca_1(\overline{\mathcal{M}}_k )^2
= \frac{k^2}{2} \zeta_K(-1)  \left(
    \frac{\zeta_K'(-1)}{ \zeta_K(-1)} + \frac{\zeta'(-1)}{\zeta(-1)} +
    \frac{3}{2} + \frac{1}{2} \log(D) \right)
\cdot E_2 (\tau,0) ,
\end{align}
where $E_2 (\tau,0)$ denotes the Eisenstein series defined in \eqref{eis}. 
In particular, the arithmetic self intersection number of
$\overline{\mathcal{M}}_k$ is  given by:
\begin{align}\label{eq:asi}
\overline{\mathcal{M}}_k^3
& = -k^3 
\zeta_K(-1)  \left(
    \frac{\zeta_K'(-1)}{ \zeta_K(-1)} + \frac{\zeta'(-1)}{\zeta(-1)} +
    \frac{3}{2} + \frac{1}{2} \log(D) \right).
\end{align}
\end{theorem}

Let us briefly indicate how Theorem \ref{thm:B} can be proved (see \cite{BBK} Theorem 6.4 for details).
For simplicity we assume that $M_2^+(D,\chi_D)=\C E_2 (\tau,0)$ (that is $D=5$, $13$, or $17$). In this case, regarding Theorem  \ref{thm:ags}, we only have to determine the constant term of $\widehat{A}(\tau) \cdot 
\cca_1(\overline{\mathcal{M}}_k )^2$, that is, essentially the arithmetic self intersection number of $\overline{\mathcal{M}}_k$.
The hypothesis on $M_2^+(D,\chi_D)$ implies in particular that $G_m$ is the logarithm of the Petersson metric of a holomorphic Borcherds product with divisor $T(m)$ for any $m$.

Let $p$ be any prime that is split in $\OK$ (that is $\chi_D(p)=1$). 
It can be shown that there exist infinitely many $m_2$ and infinitely many $m_3$, such that $\chi_D(m_j)=1$, $T(m_2)$ is disjoint to the boundary, and such that all possible intersections of $T(p)$, $T(m_2)$, $T(m_3)$ on $X_K$ are proper.

Let $F_1$, $F_2$, $F_3$ be the Borcherds products on $X_K$ with divisors $T(p)$, $T(m_2)$, $T(m_3)$, respectively. By the bilinearity of the arithmetic intersection pairing we may assume that these Borcherds products are integral, i.e., have trivial multiplier system and integral rational Fourier coefficients. 
We may further assume that they all have the same (sufficiently divisible) weight $k$.
Thus they define sections of $\calM_k$.  
The definition of the arithmetic self intersection number then says:
\begin{align}
\nonumber
\overline{\mathcal{M}}_k^3
&= \dega \left(  h_* \big(\dv(F_1)  \cdot \dv(F_2) \cdot \dv(F_3)\big) \right)\\
&\phantom{=}{} + \frac{1}{(2\pi i)^2}
\int\limits_{\widetilde{\calX}_K(\C)} (-\log\|F_1\|_\Pet)  *(-\log\|F_2\|_\Pet) *(-\log\|F_3\|_\Pet) .\label{eq:star-t(p)} 
\end{align}
Here the integral is over the star product of the Green functions corresponding to the sections $F_j$ of $\overline{\mathcal{M}}_k$. It describes the intersection at the Archimedian place. Moreover, $h:\widetilde{\calX}_K\to \Spec(\Z)$ denotes the structure morphism. The first summand is the intersection at the finite places.

Using growth estimates for certain boundary terms, one finds that the integral in \eqref{eq:star-t(p)} is equal to
\[
k^2 \int\limits_{\widetilde{X}_K} 
G_p \, \Omega^2  
+ k \int\limits_{T(p)'} -\log\|F_2\|_{\Pet} \,\Omega
+ \int\limits_{T(p)' \cap \dv(F_2)'} -\log\|F_3\|_{\Pet},
\]
where $T(p)'$ denotes the strict transform of the divisor $T(p)$ in $\widetilde{X}_K$ (\cite{BBK} Theorem 3.13).
The integral of $G_p$ was computed in \eqref{greenint}.

There is a birational morphism $\varphi$ from the modular curve $X_0(p)$ onto $T(p)$ (which extends to integral models over $\Z[1/p]$). This fact can be used to interpret the sum of the latter two integrals as a star product on the modular curve $X_0(p)$, where it can be evaluated by means of the results of \cite{Kue2} or \cite{Bost}. 

The finite intersection can also be reduced to a finite intersection on the minimal regular model of $X_0(p)$ by applying the projection formula for the morphism $\varphi$. 

It turns out that the finite contribution and the Archimedian contribution fit together rather nicely and yield the desired result up to contributions from the fiber above $p$. But now we can vary $p$, that is, take different Borcherds products for the $F_j$, to get the precise formula for $\overline{\calM}_k^3$.

In the general case, one can argue similarly, since it can be proved that
$\ch_{HZ}^1(\widetilde{X}_K)_\Q$ is already generated by Hirzebruch-Zagier divisors $T(p)$ of prime discriminant $p$ (with $\chi_D(p)=1$), see \cite{BBK} section 4.2.

\medskip

Formula \eqref{eq:asi} provides evidence for a conjecture of Kramer, based on results obtained in \cite{kramer} and \cite{Kue1}, saying that
the arithmetic volume of an arithmetic variety as $\widetilde{\calX}_K$ is essentially the derivative of the zeta value
for the geometric volume of $\widetilde{\calX}_K(\C)$. In the same way, it provides further evidence for the conjecture of Kudla on the constant term of
the derivative of certain Eisenstein series \cite{Ku:Bourbaki}, 
\cite{Ku:ICM}, \cite{Ku:MSRI}, and the conjecture of Maillot and Roessler
on special values of logarithmic derivatives of Artin  $L$-functions 
\cite{MaRo}.

It would be very interesting to find a more conceptual explanation for the fact that the geometric intersection $A(\tau)\cdot  \cc_1(\calM_k(\C))$ is proportional to $\widehat{A}(\tau) \cdot 
\cca_1(\overline{\mathcal{M}}_k )^2$.

We may apply Theorem \ref{thm:B} and $\eqref{greenint}$ to compute the Faltings height of $\calT(m)$ with respect to $\overline{\mathcal{M}}_k$ (as defined in \cite{BKK} and \cite{BBK} section 1). 
We find:

\begin{theorem}\label{thm:Tmheight} 
Recall Assumption \ref{ass:model}.  If $T(m)$ is a Hirzebruch-Zagier
divisor which is disjoint to the boundary of $X_K$, then the Faltings
height of its model $\calT(m) \in \z^1(\widetilde{\calX}_K)$ is given
by
\begin{align*}
\fh_{\overline{\mathcal{M}}_k}
(\mathcal{T}(m))&= - 2k^2 \vol_\Omega(T(m))
 \left( \frac{\zeta'(-1)}{\zeta(-1)} + \frac{1}{2}
+\frac{1}{2}\frac{\sigma_m'(-1)}{\sigma_m(-1)} 
\right).
\end{align*}
Here $\sigma_m(s)$ is the generalized divisor sum defined in $\eqref{def:sig}$.
\end{theorem}

We conclude by noticing that Assumption \ref{ass:model} can be avoided in the above theorems by introducing a level structure to rigidify the moduli problem. For instance, in \cite{BBK} the full level $N$-structure is used (where $N$ is an arbitrary integer $\geq 3$). Then the moduli problem is represented by an arithmetic variety over $\Z[\zeta_N, 1/N]$. 
However, since $N$ is inverted in the base, one only gets arithmetic intersection numbers in $\R_N = \R \big/ \big\langle \sum_{p|N} \Q \cdot\log(p) \big\rangle$.

%
%

\end{document}